\documentclass[CJK,11pt,twoside, openany]{ctexart}
\usepackage{latexsym}
\usepackage{amsthm}
\usepackage{bbm}
 \usepackage{mathrsfs}
 \usepackage{fancyhdr}
\usepackage{mathrsfs}
\usepackage{cancel}
\usepackage{color}
\usepackage{multirow}
\usepackage{eso-pic}
\usepackage{dsfont}
\usepackage{amssymb}
\usepackage{CJKspace}

\def\nnb{\nonumber}
\def\ds{\displaystyle}
\def\cd{\cdot}

\def\all{  \, \forall \, }

\def\disp{\displaystyle}

\newcommand{\refeq}[1]{~$(\ref{#1})$}
\newcommand{\eqref}[1]{~$(\ref{#1})$}
\newcommand{\thb}[1]{~{\rm (#1)}}
\def\endpf{\hfill$\Box$\vspace{0.4cm}}
\def\eqif{\, {\rm if}\, \,}

\def\eqin{ \, {\rm in } \, \,}
\def\eqae{ \, {\rm a.e. } \,\, }

\def\ul{\underline}
\def\ol{\overline}
\def\s{\sigma}

\def\Go{\omega}

\def\Ga{\alpha}
\def\Gb{\beta}
\def\Gg{\gamma}
\def\Gd{\delta}

\def\Gve{\varepsilon}
\def\Gz{\zeta}
\def\Gs{\sigma}
\def\Gt{\theta}
\def\Gvp{\varphi}

\def\GT{\Theta}

\def\GO{\Omega}
\def\bs{\bar s}
\def\bt{\bar t}
\def\bu{\bar u}
\def\bv{\bar v}

\def\by{\bar y}
\def\bz{\bar z}
\def\cM{{\cal M}}

\def\mcP{{\mathscr P}}

\def\mcR{{\mathscr R}}

\def\mcU{{\mathscr U}}
\def\mcV{{\mathscr V}}

\def\bGs{{\bar \s}}
\def\bpsi{{\bar \psi}}
\def\bGvp{{\bar \varphi}}

\def\bGP{{\ol{\Phi}}}
\def\hv{\hat v}
\def\hw{\hat w}

\def\hy{\hat y}

\def\hGs{\hat \Gs}
\def\tiM{\widetilde  M}
\def\tiN{\widetilde  N}

\def\tif{\tilde f}

\def\tis{\tilde s}
\def\tit{\tilde t}

\def\tiv{\tilde v}

\def\tix{\tilde x}
\def\tiy{\tilde y}

\def\qq{\qquad}
\def\q{\quad}
\newcommand{\dpp}[2]{{\partial {#1} \over \partial {#2}}}
\newcommand{\doo}[2]{{d {#1} \over d {#2}}}

\newcommand{\pri}{{\prime}}
\newcommand{\prii}{{\prime\prime}}

\def\supp{\,{\rm supp}\,}
\def\diam{\,{\rm diam}\,}

\def\pa{\partial}

\def\coh{\,{\rm \overline{\, co}}\,}

\def\diam{\,{\rm diam}\,}

\def\limsup{\mathop{\overline{\rm lim}}}
\def\liminf{\mathop{\underline{\rm lim}}}

\def\IR{\mathds{R}}
\def\defeq{\buildrel \triangle \over =}

\newcommand{\set}[1]{\left\{#1\right\}}
\newcommand{\ip}[1]{\left\langle #1\right\rangle}

\def\thebibliography#1{\center{\bf\normalsize References}\list
 {[\arabic{enumi}]}{\settowidth\labelwidth{[#1]}\leftmargin\labelwidth
 \advance\leftmargin\labelsep
 \usecounter{enumi}}
 \def\newblock{\hskip .11em plus .33em minus .07em}
 \sloppy\clubpenalty4000\widowpenalty4000
 \sfcode`\.=1000\relax}

\makeatletter
\def\cleardoublepage{\clearpage\if@twoside \ifodd\c@page\else
   \hbox{}\thispagestyle{empty}\newpage\addtocounter{page}{-1}
   \if@twocolumn\hbox{}\newpage\fi\fi\fi}
\makeatother
\usepackage{graphicx}
\newtheorem{Definition}{Definition}[section]
\newtheorem{Theorem}[Definition]{Theorem}
\newtheorem{Lemma}[Definition]{Lemma}
\newtheorem{Corollary}[Definition]{Corollary}
\newtheorem{Proposition}[Definition]{Proposition}
\newtheorem{Remark}{Remark}[section]

\begin{document}
\title{\bf Time optimal control problems for some non-smooth systems\thanks{This work was supported in part by 973 Program (No. 2011CB808002)  and  NSFC (No. 11371104).}}

\author{Hongwei Lou\footnote{School
of Mathematical Sciences, and LMNS, Fudan University, Shanghai
200433, China (Email: \texttt{hwlou@fudan.edu.cn}).},
~~Junjie Wen\footnote{School of Mathematical Sciences, Fudan University, Shanghai 200433, China (Email: \texttt{10110180036@fudan.edu.cn}).}
~~and~~
 Yashan Xu\footnote{School
of Mathematical Sciences, and LMNS, Fudan University, Shanghai
200433, China (Email: \texttt{yashanxu@fudan.edu.cn}).}
}

\date{}

\maketitle

\begin{quote}
\footnotesize {\bf Abstract.} Time optimal control problems for some non-smooth systems in general form are considered. The non-smoothness is caused by singularity. It is proved that Pontryagin's maximum principle holds for at least one optimal relaxed control.
Thus, Pontryagin's maximum principle holds when the optimal classical control is a unique optimal relaxed control. By constructing an auxiliary controlled system which admits the original optimal classical control as its unique optimal relaxed control,
one get a chance to get Pontryagin's maximum principle for the original optimal classical control.
Existence results are also considered.


\end{quote}

\normalsize

\def\theequation{1.\arabic{equation}}
\setcounter{equation}{0} 
\setcounter{Definition}{0} \setcounter{Remark}{0}
 \vspace{6mm}

 \section{Introduction.}

Consider the following non-smooth controlled system:
\begin{equation}\label{E101}
   \left\{\begin{array}{ll}\ds  \doo {y(t)} t=f(t,y(t),u(t)),&　t>0,\\
   y(0)=y_0, & \end{array}\right.
\end{equation}
where the state function $y(\cd)$ is a vector-valued function, while the control function $u(\cd)$ values in some metric space $U$.

In this paper, the non-smoothness means that $f$ may singular at target set. More precisely, we will consider cases of what $f$ is locally bounded while $f_y$ may unbounded near the target set.

We say $y(\cd)\in C([0,T);\IR^n)$ is a solution of \refeq{E101} on $[0,T)$ if $y(\cd)$ satisfies
\begin{equation}\label{F102}
y(t)=y_0+\int^t_0f(s,y(s),u(s))\, ds, \qq t\in [0,T).
\end{equation}

Let the target set $Q$ be a closed subset of $\IR^n$. Denote
\begin{equation}\label{F108}
\begin{array}{l}\ds
\mcU=\set{u(\cd):[0,+\infty)\to U\Big|\, u(\cd)\q\mbox{is measurable}\,},\\
\mcP=\set{(T,y(\cd),u(\cd))\in (0,+\infty)\times C([0,T);\IR^n)\times \mcU\Big| \mbox{\refeq{E101} \,holds on }\, [0,T)},\\
\ds \mcP_{ad}=\set{(T,y(\cd),u(\cd))\in \mcP\Big| \liminf_{s\to T^-}d(y(s),Q)=0},\\
\ds \mcU_{ad}=\set{u(\cd)\in \mcU\Big| (T,y(\cd),u(\cd))\in \mcP_{ad}}.
\end{array}
\end{equation}
The element in $\mcP$, $\mcP_{ad}$ and $\mcU_{ad}$ is called
feasible classical triple , admissible triple and admissible
control, respectively.

 The time optimal control problem is

\textbf{Problem (T)}: Find $(\bt,\by(\cd),\bu(\cd))\in \mcP_{ad}$  such that
\begin{equation}\label{F109}
\bt=\inf_{(T,y(\cd),u(\cd))\in \mcP_{ad}}T.
\end{equation}

\begin{Remark}\label{R101} We can see that the description above is a little different from a general one, such as the state function belongs to $C([0,T);\IR^n)$ but not $C([0,T];\IR^n)$.  On the other hand, we use $\ds \liminf_{s\to T^-}d(y(s),Q)=0$ instead of $y(T)\in Q$ to describe that the state ``reaches" the target. That is to avoid unnecessary inconvenience. Under assumptions of this paper,
it is possible that $\ds \lim_{s\to T^-} y(s)$ does not exist and then
$y(\cd)$ becomes very complicated in approaching $Q$.  Naturally, when $\ds \lim_{s\to T^-} y(s)$ exist, $\ds \lim_{s\to T^-}d(y(s),Q)=0$ becomes $y(T)\in Q$.
\end{Remark}

This paper is stimulated by Lin and Wang \cite{LinWang} for optimal blowup time problems, and Lin \cite{Lin1} for optimal quenching time problems.
In \cite{Lin1}, $f$ itself is unbounded near the target. While in \cite{LinWang}, the ``target set" is $\set{\infty}$, and $f$ is also unbounded near the target. We observed that both optimal quenching time problems and optimal blowup time problems can be transformed to time optimal control problems for some non-smooth systems (see Problem (T)) by some suitable state transform.

As far as we know, the concept of quenching appeared in  H. Kawarada  \cite{Kawarada} at the earliest, which studied the initial boundary value of the diffusion equation
\begin{eqnarray}\label{pre4}
   \left\{\begin{array}{l}
   \disp u_t= u_{xx}+\frac{1}{1-u}, \q t>0, \, x\in(0, l),\\
   u(t, 0)=u(t, l)=0, \q t>0,\\
   u(0, x)=0, \q x\in(0, l).
   \end{array}\right.
\end{eqnarray}
The author give the definition of quenching: namely, $u$ can get to 1 in finite time. The author also get some sufficient conditions of the quenching behavior.

The quenching phenomenon can be observed in the transient current of polarized ion conductor . The quenching of solution means that the derivative of the solution goes to infinity in finite time while it keep bounded itself. Many people are interested in quenching behavior of partial differential equations. See \cite{Chan1}, \cite{Chan2}, \cite{Guo} for the reference.

Blowup is another concept related to quenching which means the solution is unbounded in finite time. The difference between them is that solution keep bounded at quenching time while explode at the blowup time. We refer the readers the following works on blowup behaviors:
 \cite{Bandle}, \cite{Escobedo}, \cite{Glassey}, \cite{Guo}, \cite{Yordanov} and \cite{Zhang}, for examples.

Our aim is to study the existence theory and the necessary conditions for solutions of Problem (T).
The non-smoothness bring much more difficulties than what we though at the beginning. In this paper, we get some existence results but only establish Pontryagin's maximum principle for one of relaxed optimal controls for general systems. Nevertheless, we think that for many special systems, by using such a result and the monotonicity of controlled systems, we can finally establish Pontryagin's maximum principle for any optimal control of Problem (T).

Now, let us state the problem more precisely. We make our
basic assumptions as follows:

(S1) Set $U\subset \IR^m$ is a non-empty bounded closed set;

(S2) Set $Q$ is a non-empty convex closed set in $\IR^n$, $y_0\not\in
Q$;

(S3) Function $f(t,y,u)$ is measurable in $t\in [0,+\infty)$ and continuous in $(y, u)\in (\IR^n\setminus
Q)\times U$. Moreover, for any $E\subset\subset\times
(\IR^n\setminus Q)$, there exists a constant $L_E>0$ and uniform
modulus of continuity $\Go_E(\cd)$ such that
\begin{equation}\label{E200}
    |f(t,y,u)|\leq L_E(|y|+1), \qq  \all (t,y,u)\in E\times U,
\end{equation}
\begin{equation}\label{E200}
    |f(t,y,u)-f(t,x,v)|\leq L_E|y-x|+\Go_E(|u-v|), \q  \all (t,x,u), (t,y,v)\in E\times U.
\end{equation}

(S4) Function $f(t,y,u)$ is differentiable in $y\in \IR^n\setminus Q$, for any $E\subset\subset [0,+\infty)\times
(\IR^n\setminus Q)$, there exists $M_E>0$, such that
\begin{equation}\label{E203}
    |f_y(t,y,u)|\leq M_E, \qq \all (t,y,u)\in E\times U.
\end{equation}

For $f=\pmatrix{f^1 & f^2 & \ldots & f^n}^\top$, we denote
$$
f_t=\pmatrix{{\pa f^1\over \pa t} \cr {\pa f^2\over \pa t}\cr \vdots\cr {\pa f^n\over \pa t}}, \q f_y= \dpp f y= \pmatrix{
  {\pa f^1\over \pa y_1} &  {\pa f^2\over \pa y_1} & \cdots & {\pa f^n\over \pa y_1}\cr
   {\pa f^1\over \pa y_2} &  {\pa f^2\over \pa y_2} & \cdots & {\pa f^n\over \pa y_2}\cr
   \vdots & \vdots & \ddots & \vdots \cr
     {\pa f^1\over \pa y_n} &  {\pa f^2\over \pa y_n} & \cdots & {\pa f^n\over \pa y_n}}.
$$

As mentioned above, $y(\cd)$ is not always well-defined
after it reaches $Q$. Since we only concern about the first time
when it reaches $Q$, and $y(\cd)$ is uniquely determined by
\refeq{E101} before that time from (S3)---(S4), we can denote the
solution of \refeq{E101} by $y(\cd;u(\cd))$ without any
misunderstanding.

The outline of this paper is: in first section, we give a general description of the problem. In the second section, we explain the problem in detail and make some transformation. The third section will devote the relaxed control. The fourth section will give the existence of the optimal control, and the fifth section will focus on the maximum principle of the optimal control.

\def\theequation{2.\arabic{equation}}
\setcounter{equation}{0} 
\setcounter{Definition}{0} \setcounter{Remark}{0}
\section{Relaxed Problems and the Corresponding Results.}  To study the existence of optimal control, we
 introduce the relaxed control which can also help to study the necessary
  conditions of optimal controls. To deal with the non-smoothness of
the system near the target set, we need to introduce the
corresponding approximate problem. However,  approximate problems need not necessary admits optimal control even if the optimal control of the original problem exists. To overcome this difficulty, we consider relaxed controls.

\subsection{Relaxed controls.}
Now,
we recall the notion of relaxed control and state some preliminary
results about the space of relaxed controls.

We denote by
$\cM^1_+(U)$ the set of all  probability measures in
$U$, by $\mcR(U)$ ($\mcR_T(U)$ )
the set of all measurable probability
measure-valued functions on $[0,+\infty)$ ($[0,T]$), that is, $\Gs(\cd)\in \mcR(U) $  ($\mcR_T(U)$ ) if and only if
$$
 \Gs(t)\in  \cM^1_+(U), \qq \eqae t\in [0,+\infty) \,([0,T]),
$$
and
$$
t\mapsto \int_U h(v)\Gs(t)(dv) \mbox{ is measurable,}\qq \all
h\in C(U),
$$
where $C(U)$ denotes the space of all continuous functions on $U$. Let
$C(U)^*$ and $L^1([0,T];C(U))^*$ be the dual spaces of $C(U)$ and $L^1([0,T];C(U))$ with weak star topology, respectively. We
regard  $\cM^1_+(U)$ and  $\mcR_T(U) $ as
subspaces of $C(U)^*$ and $L^1([0,T];C(U))^*$, respectively,
by setting
\begin{equation}
\label{E29} \Gt(h)\defeq \int_U h(v)\Gt(dv), \qq \all \Gt(\cd)\in
        \cM^1_+(U), \qq\all h\in C(U)
\end{equation}
and
\begin{eqnarray}
\nnb
& & \Gs(g)\defeq\int^T_0 dt\int_U g(t,v)\Gs(t)(dv), \\
& & \label{E210a}  \mbox{\hspace{2.5cm}}
     \all  \Gs(\cd)\in  \mcR(U),\, g\in  L^1([0,T];C(U)).
\end{eqnarray}
We see that   \refeq{E210a} is well-defined by Theorem {\rm IV.}1.6, (p. 266) in  \cite{W3}. We say that
$$
\Gs_k(\cd)\to \Gs(\cd), \qq\eqin   \mcR_T(U)
$$
if
\begin{eqnarray*}
& &\int^T_0ds\int_U h(s,v)\Gs_k(s)(dv)\to
\int^1_0ds\int_U h(s,v)\Gs(s)(dv),\\
& &    \qq\qq\qq\qq\all  h\in  L^1([0,T];C(U))
\end{eqnarray*}
while say that
$$
\Gs_k(\cd)\to \Gs(\cd), \qq\eqin   \mcR(U)
$$
if
$$
\Gs_k(\cd)\to \Gs(\cd), \qq\eqin   \mcR_T(U), \q\all T>0.
$$
In fact, we only concern $\mcR_T(U)$ for $T$  large enough, but not
$\mcR(U)$ itself. We introduce $\mcR(U)$ just for the convenience
as the optimal time is unknown for us.

The following lemma is an important property of relaxed controls.
\begin{Lemma}\label{T21} Assume $T>0$ and $U$ be a compact metric space.
Then $\mcR_T(U)$ is convex and sequentially compact.
\end{Lemma}
For a proof of the above lemma, see Warga \cite{W3}, Theorem
{\rm IV}.2.1, p. 272.

\subsection{Existence of optimal relaxed controls.}

First, we study the existence of the optimal relaxed controls. Consider the following relaxed controlled system corresponding to
\refeq{E101}:
\begin{equation}\label{E204a}
   \left\{\begin{array}{ll}\ds  \doo {y(t)} t=\int_U  f(t,y(t),v)\Gs(t)(dv),&　t>0,\\
   y(0)=y_0.  & \end{array}\right.
\end{equation}
Denote
\def\mcRP{\mcR\!\mcP}
\begin{equation}\label{E204A1}
\begin{array}{l}\ds
\mcRP=\set{(w,y(\cd),\Gs(\cd))\in (0,+\infty)\times C([0,w);\IR^n)\times \mcR(U)\Big| \mbox{\refeq{E204a} \,holds on }\, [0,w)},\\
\ds \mcRP_{ad}=\set{(w,y(\cd),\Gs(\cd))\in \mcRP\Big| \liminf_{t\to w^-}d(y(t),Q)=0},\\
\ds \mcR_{ad}=\set{\Gs(\cd)\in \mcR(U)\Big| (w,y(\cd),\Gs(\cd))\in \mcRP_{ad}}.
\end{array}
\end{equation}
Sets $\mcRP$, $\mcRP_{ad}$ and $\mcR_{ad}$ are called the set
of feasible relaxed triples , the set of admissible relaxed triples
and the set of admissible relaxed controls.

The optimal relaxed control problem is

\textbf{Problem ($\mcR$).} Find $(w^*,y^*(\cd),\Gs^*(\cd))\in \mcRP_{ad}$ such that
\begin{equation}\label{206b}
w^*=\inf_{(T,y(\cd),\Gs(\cd))\in\mcRP_{ad}} w.
\end{equation}
Any triple $(w^*,y^*(\cd),\Gs^*(\cd))\in\mcRP_{ad}$ satisfying \refeq{206b} is called an optimal relaxed triple of Problem ($\mcR$)/Problem (T). If $(w^*,y^*(\cd),\Gd_{u^*(\cd)})$ is
an optimal relaxed triple of Problem ($\mcR$) for some $u^*(\cd)\in \mcU$, then we call $(w^*,y^*(\cd),u^*(\cd))$ a classical optimal triple of Problem ($\mcR$)/Problem (T).

We give the following lemma, which concerns the continuous dependence of the solutions of \refeq{E204a} with respect to the
relaxed controls in the meaning of convergence in $\mcR(U)$.
\begin{Lemma}\label{T201} Assume that \thb{S3}---\thb{S4} hold,
$\hw>0$, and $\Gs_k(\cd), \hGs(\cd)\in \mcR(U)$ such that
\begin{equation}\label{E204bb}
 \Gs_k(\cd) \to  \hGs(\cd), \qq\eqin  \mcR(U),
\end{equation}
the solution $\hy(\cd)$ of \refeq{E204a} corresponding to
$\hGs(\cd)$ exists on $[0,\hw]$, $(\hw,\hy(\cd),\hGs(\cd))\not\in
\mcRP_{ad}$. Then there exists $\Gd>0$ and $N>0$ such that the
solution $y_k(\cd)$ of \refeq{E204a} corresponding to $ \Gs_k(\cd)$
exists on $[0,\hw+\Gd]$. Moreover, $y_k(\cd)$ is uniformly bounded and equicontinuous with
respect to $k\geq N$.
\end{Lemma}
\proof Since $(\hw,\hy(\cd),\hGs(\cd))\not\in \mcRP_{ad}$, we
have
$$
\Gve_0\equiv \min_{t\in [0,\hw]}d(\hy(t),Q)>0.
$$
Denote
$$
R=\max_{t\in [0,\hw]}|\hy(t)|+1, \q \GO=\set{x\in \IR^n\Big| |x|< R, \, d(x,Q)> {\Gve_0\over 2}}, \q  E=[0,\hw+1]\times \GO.
$$
Let $\ds \Gd={\min(\Gve_0, 1)\over 12  e^{ L_E(\hw+1)}}$.  Taking
an integer $\ell>{L_E(R+2)\over \Gd}$. By (S3), we get
\begin{equation}\label{E205bb}
 \sup_{0\leq j\leq \ell-1\atop u\in U}\int^{(j+1)\hw \over \ell}_{j\hw \over \ell} | f(t, \hy(t),u)|\, dt \leq \Gd.
\end{equation}
On the other hand, it follows from \refeq{E204bb} that there exists $N>0$, such
that when $k\geq N$,
\begin{equation}\label{E206bb}
\Big|\int^{j\hw \over \ell}_0dt\int_U f(t, \hy(t),u)\,\big(\Gs_k(t)-\hGs(t)\big)(du)\Big| \leq \Gd,\q \all j=1,2,\ldots,\ell-1.
\end{equation}
Combining the above with \refeq{E205bb}, we get
\begin{equation}\label{E207bb}
\Big|\int^t_0ds\int_U f(s, \hy(s),u)\,\big(\Gs_k(s)-\hGs(s)\big)(du)\Big| \leq 3\Gd,\q \all t\in [0,1], \q k\geq N.
\end{equation}
By (S3), we see that the solution of \refeq{E204a} exists before
it reaches the target set $Q$. We assert that for $k\geq N$, the
corresponding $y(\cd)=y(\cd;\Gs_k(\cd))$ exists on $[0,\hw
+\Gd]$ and satisfies
\begin{equation}\label{E207cc}
y(t) \in \GO, \qq t\in [0,\hw+\Gd], \q k\geq N.
\end{equation}
Otherwise, there exists $S\in (0,\hw+\Gd)$ such that $y(\cd)=
y(\cd;\Gs_k(\cd))$ exists on $[0,S]$ and satisfies
\begin{equation}\label{E208}
y(S)\in \partial \GO, \q y(t)\in \GO, \qq\all t\in [0,S).
\end{equation}
We have
\begin{eqnarray*}
&& |y(t)-\hy(t)|\\
&=& \Big|\int^t_0ds\int_U  f(s,y(s),u)\Gs_k(s)(du)-\int^t_0ds\int_U f(s,\hy(s),u)\hGs_k(s)(du)\Big|\\
&\leq &  \Big|\int^t_0ds\int_U  \Big(f(s,y(s),u)-f(s,\hy(s),u)\Big)\Gs_k(s)(du)\Big|\\
&&  +\Big|\int^t_0ds\int_U  f(s,\hy(s),u)\Big(\Gs_k(s)-\hGs(s)\Big)(du)\Big|\\
&\leq &  L_E\int^t_0 |y(s)-\hy(s)| \, ds+3 \Gd,\qq\all t\in  [0,S].
\end{eqnarray*}
Then it follows from Grownwall's inequality that
\begin{equation}\label{E216A}
|y(t)-\hy(t)|\leq  3 \Gd e^{L_E t}\leq  {1\over 4}\min(\Gve_0, 1), \qq\all t\in [0,S].
\end{equation}
Particularly, we have
$$
|y(S)|\leq |\hy(S)|+{1\over 4}<R, \q d(y(S),Q)\geq {3\Gve_0\over 4},
$$
which contradicts to \refeq{E208}. Therefore,
$y(\cd)=y(\cd;\Gs_k(\cd))$ exists on $[0,\hw+\Gd]$ . Furthermore, similar to the proof of \refeq{E216A}, it can be proved that
$y(\cd;\Gs_k(\cd))$
 is uniformly bounded respect to $k\geq N$ . Then we get also the equicontinuity. We reach the conclusion.
\endpf

From the above lemma we conclude that under assumptions
(S3)---(S4) , if $y_0\not\in Q$, then the optimal time of Problem
$(\mcR)$ is positive.

\begin{Lemma}\label{T202} Assume that \thb{S3}---\thb{S4} hold. If there exists $(w_k,y_k(\cd),\Gs_k(\cd))\in
\mcRP$ such that
\begin{equation} \label{E218}
 \lim_{k\to +\infty}w_k = w, \q \lim_{k\to +\infty}\lim_{t\to w_k^-}d(y_k(t), Q)= 0.
\end{equation}
Then the optimal time of Problem $(\mcR)$ is not bigger than $w$.
\end{Lemma}

\proof  By Lemma \ref{T21}, we can suppose that
\begin{equation}\label{E219}
 \Gs_k(\cd) \to  \Gs(\cd), \qq\eqin  \mcR(U).
\end{equation}
Then, it follows easily from Lemma \ref{T201} that there must exist $S\in (0,w]$ such that
$y(\cd;\Gs(\cd))$ exists on $[0,S)$ and
$$
 \liminf_{t\to S^-}d(y(t;\Gs(\cd)),Q)=0.
$$
That is, $(S,\hy(\cd),\hGs(\cd))\in \mcRP_{ad}$. Therefore the optimal
time of Problem $(\mcR)$ is not bigger than $S\leq w$.
\endpf

\begin{Remark}\label{R201}
The lemma above says that if the optimal time of Problem $(\mcR)$ is
$t^*$, then for any $S\in (0,t^*)$, and $\Gs(\cd)\in \mcR$, the
solution of \refeq{E204a} exists on $[0,S]$. Moreover, the distance
between
$$
P_S\equiv \set{y(\cd;\Gs(\cd))\Big|_{[0,S]}\,\Big|\Gs(\cd)\in \mcR}
$$
and $Q$ is positive, that is, $\ds \inf_{t\in [0,S]\atop y(\cd )\in
P_S}d(y(t),Q)$ is positive. Consequently, $P_S$ is uniformly bounded and
equicontinuous on $[0,S]$.
\end{Remark}
\begin{Remark}\label{R202}
Lemma \ref{T201} shows that under \thb{S3}---\thb{S4}, $\mcRP_{ad}\ne \emptyset$ is equivalent to
$$
\inf_{(w,y(\cd),\Gs(\cd))\in \mcRP}\liminf_{t\to w^-}d(y(t),Q)=0.
$$
\end{Remark}

The following theorem gives the existence of optimal  relaxed  triples under relatively weak assumptions.
\begin{Theorem}\label{T203}  Assume that \thb{S1}---\thb{S3} hold. Moreover, suppose that $y_0\not\in Q$ and
 $\mcRP_{ad}\ne \emptyset$. Then Problem \thb{$\mcR$} admits at least one optimal relaxed  triple.
\end{Theorem}
\proof The proof is quite standard. Let $(w_k,y_k(\cd),\Gs_k(\cd))\in \mcRP_{ad}$ be a minimizing sequence, i.e.,
$$
\lim_{k\to \infty}w_k=w^*\equiv \inf_{(w,y(\cd),\Gs(\cd))\in\mcRP_{ad}} w.
$$
By Lemma \ref{T21}, we can suppose that:
\begin{eqnarray*}
\Gs_k(\cd)&\to &\Gs^*(\cd), \qq \eqin \mcR_{w^*}(U).
\end{eqnarray*}
Then Remark \ref{R201} shows that $y_k(\cd)\equiv y(\cd;\Gs_k(\cd))$ is uniformly
bounded and  equicontinuous on $[0,w]$ for any $w\in (0,w^*)$.  Thus by Arzel\'a-Ascoli's theorem,
we can suppose that $y_k(\cd)$
converges uniformly in $C([0,w];\IR^n)$  to some $y^*(\cd)\in C([0,w^*);\IR^n)$. By (S4), there is a constant $C_w>0$  such that
\begin{eqnarray*}
&& \Big|\int^t_0 ds\int_U   f(s, y_k(s), v)\Gs_k(s)(dv)-\int^t_0 ds\int_U f(s, y^*(s), v)\Gs^*(s)(dv)\Big|\\
&\leq &\Big|\int^t_0 ds\int_U f(s, y^*(s), v)\Big(\Gs_k^*(s)-\Gs^*(s)\Big)(dv)\Big|\\
&& +C_w\int^t_0\int_U  \Big| y_k(s)-y^*(s)\Big|\, ds, \q\all t\in [0,w].
\end{eqnarray*}
Therefore,
$$
\lim_{k\to +\infty}\Big|\int^t_0 ds\int_U   f(s, y_k(s), v)\Gs_k(s)(dv)-\int^t_0 ds\int_U f(s, y^*(s), v)\Gs^*(s)(dv)\Big|=0, \q\all t\in [0,w^*).
$$
Consequently
$$
y^*(t)=y_0+ \int^t_0 ds\int_U f(s, y^*(s), v)\Gs^*(s)(dv), \q\all t\in [0,w^*).
$$
That is, $(w^*,y^*(\cd),\Gs^*(\cd))\in \mcRP$. On the other hand, using Lemma \ref{T201}, we can get that $(w^*,y^*(\cd),\Gs^*(\cd))\in \mcRP_{ad}$. Consequently, $(w^*,y^*(\cd),\Gs^*(\cd))$ is an optimal relaxed triple.
\endpf

\subsection{Maximum Principle for Optimal Relaxed Triple.}

To yield the maximum principle for optimal relaxed triple, we introduce approximate problems.

For any $\Ga>0$, denote
$$
Q_\Ga=\set{q\in \IR^n|d(q,Q)\leq  \Ga }.
$$
Then  $Q_\Ga$ is closed and convex.  Moreover, denote $\Ga_0=d(y_0,Q)$, then when $\Ga\in (0,\Ga_0)$, $y_0\not\in Q_\Ga$.
Let
$$
\mcRP^\Ga_{ad}=\set{(w,y(\cd),\Gs(\cd))\in\mcRP | y(w)\in Q_\Ga}.
$$

We introduce the following approximate problems.

\textbf{Problem ($\mcR_\Ga$).}  Find $(w^{\Ga,*},y^{\Ga,*}(\cd),\Gs^{\Ga,*}(\cd))\in\mcRP^\Ga_{ad}$ such that
\begin{eqnarray}\label{206bc}
\nnb && w^{\Ga,*}=\inf_{(w,y(\cd),\Gs(\cd))\in\mcRP^\Ga_{ad}}  w.
\end{eqnarray}

We have the following result.
\begin{Lemma}\label{T205}  Assume that \thb{S1}---\thb{S4}  hold and $\mcRP_{ad}\ne \emptyset$.  Let  $\Ga\in (0,\Ga_0)$. Then Problem \thb{$\mcR_\Ga$} admits at least one optimal triple  $(w^{\Ga,*},y^{\Ga,*}(\cd),\Gs^{\Ga,*}(\cd))$.

Moreover, there exists a nontrivial solution of
\begin{equation}\label{E227a}
   \doo {\psi^{\Ga,*}(t)} t= \ds -\int_U f_y (t,y^{\Ga,*}(t),v)  \psi^{\Ga,*}(t)\Gs^{\Ga,*}(t)(dv),\q    t\in [0,w^{\Ga,*}),
\end{equation}
such that
\begin{equation}\label{E231a}
\ip{ \psi^{\Ga,*}(w), q-y^{\Ga,*}(w)}\geq 0, \qq\all q\in Q_\Ga
\end{equation}
and
\begin{eqnarray}\label{E231}
  && \supp \Gs^{\Ga,*}(t)\subseteq \Big\{\tiv\in U\left| \begin{array}{l}\ds \qq  H(t,y^{\Ga,*}(t),\psi^{\Ga,*}(t),\tiv)\\
\ds =\max_{v\in U }H(t,y^{\Ga,*}(t),\psi^{\Ga,*}(t),v)\end{array}\right.\Big\}, \, \eqae t\in [0,w^{\Ga,*}),
\end{eqnarray}
where
\begin{equation}\label{E232}
H (t,y,\psi,v)=   \ip{\psi, f (t,y,v)}.
\end{equation}
\end{Lemma}
In the above lemma, the existence results can be looked a corollary of Theorem \ref{T203}, while the necessary conditions can be yielded in a standard way like that for optimal classical triple in smoothness cases.

\begin{Remark}\label{R305} Let $\GO\supset\supset Q_\Ga$ and $G\in C^1(\GO;\IR^n)$. Moreover, $\ds {\pa G(x)\over \pa x}$ is not singular
for any $x\in \GO$. If
$\ds \mcRP^\Ga_{ad}$ is replaced by
$$
\mcRP^\Ga_{ad}=\set{(w,y(\cd),\Gs(\cd))\in\mcRP | G(z(w))\in Q_\Ga},
$$
then Lemma \ref{T205} holds with \refeq{E231a} being replaced by
\begin{equation}\label{E231b}
\ip{ \Big({\pa G\over \pa x}(y^{\Ga,*}(w))\Big)^{-1}\psi^{\Ga,*}(w), q-G(y^{\Ga,*}(w))}\geq 0, \qq\all q\in Q_\Ga.
\end{equation}
\end{Remark}

The following theorem is related to the necessary conditions for optimal relaxed controls. To get the transversality condition, we set the following assumptions.

(S3$^\pri$)  Function $f(t,y,u)$ is measurable in $t\in [0,+\infty)$ and continuous in $(y, u)\in (\IR^n\setminus
Q)\times U$. Moreover, there exists a constant $L>0$ such that
\begin{equation}\label{E200}
    |f(t,y,u)|\leq L(|y|+1), \qq  \all (t,y,u)\in [0,+\infty)\times (\IR^n\setminus Q)\times U.
\end{equation}

(CE1)   There exists   $\Gd_1>0$ such that
$\ds  f_y(t,y,u)+f_y(t,y,u)^\top $
is positive simi-definite for any $t\geq 0$, $u\in U$ and $d(y,Q)\leq \Gd_1$.

(CE2) For any $\Gve>0$, there exists $C_\Gve>0$,   $\Gd_\Gve>0$, $m_\Gve\geq 1$  and  $E_1,E_2,\ldots, E_{m_\Gve}\subset S^{n-1}$ such that
\begin{equation}\label{EEEcde1}
 \diam E_k\equiv \sup_{x,y\in E_k}|x-y| \leq \Gve, \qq k=1,2,\ldots, m_\Gve,
 \end{equation}
 \begin{equation}\label{EEEcde2}
d(E_k,E_j)\geq \Gve, \qq 1\leq k<j\leq m_\Gve,
 \end{equation}
and
\begin{eqnarray}\label{EEEaaaa}
\nnb  && |f_y(t,y,u)\xi|\leq C_\Gve\ip{f_y(t,y,u)\xi, \xi}+C_\Gve,  \\
&& \qq\qq\q\all \xi\in S^{n-1}\setminus \bigcup^{m_\Gve}_{k=1} E_k, \, t\in [0,{1\over \Gve}],\, d(y,Q)\leq \Gd_\Gve,\, |y|\leq {1\over \Gve}.
\end{eqnarray}

\begin{Theorem}\label{T209}
Let $y_0\not\in Q$.  Assume that \thb{S1}---\thb{S4} hold. Moreover, assume that $\mcRP_{ad}\ne \emptyset$.
Then Problem $(\mcR)$ admits an optimal triple
$(w^*,y^*(\cd),\Gs^*(\cd))$ which satisfies
\begin{eqnarray}\label{E231nbb}
 && \supp \Gs^*(t)\subseteq \Big\{u\in U\left| \begin{array}{l}\ds \qq \ip{\psi^*(t), f(t,y^*(t),u)}\\
\ds =\max_{v\in U } \ip{\psi^*(t), f(t,y^*(t),v)}\end{array}\right.\Big\},\q\eqae t\in [0,w^*)
\end{eqnarray}
with $\psi^*(\cd)$ being  a nontrivial solution of
\begin{equation}\label{E227ab}
   \doo {\psi^*(t)} t = -\int_U f_y(t,y^*(t),u)  \psi^*(t)\Gs^*(t)(du),\q    t\in [0,w^*).
\end{equation}
Moreover, if $(S3^\pri)$ and \thb{CE1}---\thb{CE2} hold, then $\ds y^*(w^*)=\lim_{t\to w^{*-}}y^*(t)$ and $\ds \psi^*(w^*)=\lim_{t\to w^{*-}}\psi^*(t)$
exist,
\begin{equation}\label{E357}
\ip{ \psi^*(w^*), q-y^*(w^*)}\geq 0, \qq\all q\in Q.
\end{equation}
\end{Theorem}
\proof Let $\Ga_0=d(y_0,Q)$. Choose $\Ga\in (0,\Ga_0)$. By Lemma \ref{T205},  Problem \thb{$\mcR_\Ga$} admits at least one optimal triple  $(w^{\Ga,*},z^{\Ga,*}(\cd),\Gs^{\Ga,*}(\cd))$
satisfying \refeq{E227a}---\refeq{E231}. Moreover, we can set $|\psi^{\Ga,*}(0)|=1$.

Similar to that in the proof of Theorem \ref{T203} , there exists a sequence $\Ga_k\to 0^+$ such that
\begin{equation}\label{E248}
\left\{\begin{array}{ll}\ds
w^{\Ga_k,*}\to w^*, & \eqin \IR, \\
 \ds  y^{\Ga_k,*}(\cd)\to y^*(\cd), & \eqin C([0,w];\IR^n),\q\all w\in (0,w^*), \\
 \ds \psi^{\Ga_k,*} (\cd)\to \psi^*(\cd), & \eqin C([0,w];\IR^n), \q\all w\in (0,w^*),\\
\ds \Gs^{\Ga_k,*}(\cd)\to \Gs^*(\cd), & \eqin  \mcR(U). \end{array} \right.
\end{equation}
We can get easily that $\ds (w^*,y^*(\cd),\Gs^*(\cd))\in\mcRP_{ad}$. For any $\ds (w,y(\cd),\Gs(\cd))\in \mcRP_{ad}$, noting that $\mcRP_{ad}\subseteq\mcRP^\Ga_{ad}$  for any $\Ga\in (0,\Ga_0)$, we have
$$
w^{\Ga_k,*}\leq w,\qq \all k\geq 1.
$$
Let $k\to +\infty$, we get $w^*\leq w$.  That is, $(w^*,y^*(\cd),\Gs^*(\cd))$ is an optimal relaxed triple to Problem ($\mcR$).

By \refeq{E248}, we get from \refeq{E227a} that
\begin{equation}\label{E250}
   \doo {\psi^*(t)} t= \ds -\int_U    f_y(t,y^*(t),v)  \psi^*(t)\Gs^*(t)(dv),\q    t\in [0,w^*).
\end{equation}
Moreover, it follows from \refeq{E248} and
$$
\ds \psi^{\Ga_k,*} (\cd)\ne 0
$$
that
\begin{equation}\label{E251a}
   \psi^*(t)\ne 0,\qq \all  t\in [0,w^*).
\end{equation}
Combining \refeq{E231}  with \refeq{E248}, we get that for any $w\in (0,w^*)$ and $\Gs(\cd)\in \mcR(U)$,
\begin{eqnarray}\label{E229ccc}
\int^w_0dt\int_U \ip{\psi^*(t),  f(t,y^*(t),v)}(\Gs^*(t)-\Gs(t))(dv)\geq 0.
\end{eqnarray}
Then
\begin{eqnarray}\label{E231bb}
 && \supp \Gs^*(t)\subseteq \Big\{\tiv\in U\left| \begin{array}{l} \ds \qq H(t,y^*(t),\psi^*(t),\tiv)\\
  \ds =\max_{v\in U }H(t,y^*(t),\psi^*(t),v)\end{array}\right. \Big\}, \q \eqae t\in [0,w^*).
\end{eqnarray}

Now, let (S3$^\pri$) and \thb{CE1}---\thb{CE2}  hold. We want
yield the transversality condition \refeq{E357}.

Define
$$
z^{\Ga_k,*}(s)=y^{\Ga_k,*}(w^{\Ga_k,*}s), \qq s\in [0,1],
$$
$$
\Gvp^{\Ga_k,*}(s)=\psi^{\Ga_k,*}(w^{\Ga_k,*}s), \qq s\in [0,1].
$$
$$
\Gt^{\Ga_k,*}(s)=\Gs^{\Ga_k,*}(w^{\Ga_k,*}s), \qq s\in [0,1].
$$
We will prove that  $\Gvp^{\Ga_k,*}(\cd)$ is equicontinuous on $[0,1]$, which implies
\begin{equation}\label{Edd1}
 \Gvp^{\Ga_k,*} (\cd) \to \Gvp^*(\cd), \qq \eqin C([0,1];\IR^n).
\end{equation}
We note that it always holds that
\begin{equation}
 w^{\Ga_k,*}\leq w^*, \qq k=1,2,\ldots.
\end{equation}

By (S3$^\pri$), we can use Gronwall's inequality to get that $z^{\Ga_k,*}(\cd)$ is bounded uniformly on $[0,1]$.  That is
$$
|z^{\Ga_k,*}(s)|\leq M, \qq \all s\in [0,1], k=1,2,\ldots
$$
for some $M>0$. Then, using (S3$^\pri$)  again, we get  $\ds z^{\Ga_k,*}(\cd
)$ is equicontinuous on $[0,1]$. Thus instead of \refeq{E248}, we can easily yield that
\begin{equation}\label{E261}
   z^{\Ga_k,*}(\cd)\to z^*(\cd), \qq \eqin C([0,1];\IR^n)
\end{equation}
for some $z(\cd)\in C([0,1];\IR^n)$ with
$$
z^*(s)=y^*(w^*s), \qq s\in [0,1).
$$
Consequently $\ds y^*(\cd)$ can be well-defined and is continuous on $[0,w^*]$, i.e., $\ds y^*(w^*)=\lim_{t\to w^{*-}}y(t)$ exists. Moreover,  $\ds y^*(w^*)=\lim_{k\to +\infty}z^{\Ga_k,*}(1)=z^*(1)$.

On the other hand,  by (CE1) and noting that
$$
\ip{A\xi, \xi}\geq 0
$$
holds when $A+A^\top\geq 0$, we have
\begin{eqnarray}\label{E2}
\nnb  &&  \doo {} s \Big|\Gvp^{\Ga_k,*}(s)\Big|^2\\
\nnb &=& -2\int_U w^{\Ga_k,*} \ip{f_y(w^{\Ga,*}s,z^{\Ga_k,*}(s),v)  \Gvp^{\Ga_k,*}(s),\Gvp^{\Ga,*}(s)} \Gt^{\Ga_k,*}(s)(dv)\\
\nnb &\leq & \left\{\begin{array}{ll}\ds 0, & \eqif d(z^{\Ga_k,*}(s),Q)<\Gd_1, \\
\ds 2N_\Gve w^*\, \Big|\Gvp^{\Ga_k,*}(s)\Big|^2, & \eqif d(z^{\Ga_k,*}(s),Q)\geq \Gd_1  \end{array}\right. \\
&\leq & 2N_\Gve w^*\,  \Big|\Gvp^{\Ga_k,*}(s)\Big|^2, \qq\all s\in [0,1],
\end{eqnarray}
where
$$
N_\Gve\equiv \sup \set{|f_y(t,y,u)|\Big|  t\in [0,w^*], d(y,Q)\geq \Gd_1, |y|\leq M}<+\infty
$$
 by (S4). Combining the above with $\ds \Big|\Gvp^{\Ga_k,*}(0)\Big|=1$ we get that
$\Gvp^{\Ga_k,*}(\cd)$ is bounded uniformly on $[0,1]$.  That is
$$
\ds  \Big|\Gvp^{\Ga_k,*}(s)\Big|\leq \tiM, \qq\all s\in [0,1],k\geq 1
$$
for some $\tiM>0$. Further, we get immediately that $|\Gvp^{\Ga_k,*}(\cd)|^2$ is equicontinuous on $[0,1]$.
Therefore, we  have
\begin{equation}\label{Edd}
|\Gvp^{\Ga_k,*} (\cd)|^2\to |\Gvp^*(\cd)|^2, \qq \eqin C([0,1];\IR^n).
\end{equation}
We turn to prove \refeq{Edd1}. We mention that for any $\Ga\in (0,1)$,  $\Gvp^{\Ga_k,*}(\cd)$ is equicontinuous on $[0,1-\Ga]$ (see \refeq{E248}).

Case I: $\ds \Gvp^*(1)=0$. Then \refeq{Edd} implies \refeq{Edd1}.

Cade II:  $\ds \Gvp^*(1)\ne 0$. Then there exists $c_0>0$ such that
\begin{equation}\label{Edd2}
|\Gvp^{\Ga_k,*} (s)|\geq c_0, \qq \all s\in [0,1], k\geq 1.
\end{equation}
For any $\Gve\in (0, \max(1/w^*, 1/M))$, we get from (CE2) that
there exists $C_\Gve>0$, $\Gd_\Gve>0$, $m_\Gve\geq 1$ and  $E_1,E_2,\ldots, E_{m_\Gve}\subset S^{n-1}$ such that
\begin{equation}\label{E367}
 \diam E_k\equiv \sup_{x,y\in E_k}|x-y| \leq \Gve, \qq k=1,2,\ldots, m_\Gve,
 \end{equation}
 \begin{equation}\label{E368}
d(E_k,E_j)\geq \Gve, \qq 1\leq k<j\leq m_\Gve
\end{equation}
and
\begin{eqnarray}\label{EEE}
 \nnb && |f_y(t,y,u)\xi|\leq C_\Gve\ip{f_y(t,y,u)\xi, \xi}+C_\Gve,  \\
 && \qq\qq\qq\all \xi\in S^{n-1}\setminus \bigcup^{m_\Gve}_{k=1} E_k, \, t\in [0, w^*),\,  d(y,Q)\leq \Gd_\Gve,\,  |y|\leq M.
\end{eqnarray}
Withpout loss of genality, we soppose that $E_1, E_2, \ldots, E_{m_\Gve}$ are closed. Let
$$
\tiN=\sup \set{|f_y(t,y,u)|\Big| t\in [0, w^*), d(y,Q)\geq \Gd_\Gve, |y|\leq M }+C_\Gve
$$
and $\Gd\in (0,1)$ small enough such that
\begin{equation}\label{E369a}
2w^*\tiN\Gd+ {C_\Gve\over  c_0^2}\sup_{k\geq 1} \Big(|\Gvp^{\Ga_k,*} (1-\Gd)|^2-|\Gvp^{\Ga_k,*} (1)|^2\Big)\leq {\Gve\over 2}.
\end{equation}
Now, consider $k\geq 1$ and $1-\Gd<s_1<s_2<1$. We have
\begin{eqnarray}
\nnb && \doo {}s {\Gvp^{\Ga_k,*} (s)\over |\Gvp^{\Ga_k,*} (s)|}\\
\nnb & =&  {1\over |\Gvp^{\Ga_k,*} (s)|} \int_U w^{\Ga_k,*}\Big(I-{\Gvp^{\Ga_k,*} (s)\Gvp^{\Ga_k,*} (s)^\top\over |\Gvp^{\Ga_k,*} (s)|^2}\Big)f_y(w^{\Ga_k,*}s,z^{\Ga_k,*}(s),v)  \Gvp^{\Ga_k,*} (s) \Gt^{\Ga_k,*}(s)(dv).
\end{eqnarray}
(a) If
\begin{equation}\label{E370}
{\Gvp^{\Ga_k,*} (s)\over |\Gvp^{\Ga_k,*} (s)|}\not\in \bigcup^{m_\Gve}_{k=1}E_k, \qq\all s\in (s_1,s_2),
\end{equation}
then
\begin{eqnarray*}
\nnb && \Big|\doo {}s {\Gvp^{\Ga_k,*} (s)\over |\Gvp^{\Ga_k,*} (s)|}\Big| \\
\nnb &\leq & 2 w^{\Ga_k,*}\,\int_U \Big|f_y(w^{\Ga_k,*}s,z^{\Ga_k,*}(s),v)  {\Gvp^{\Ga_k,*} (s)\over |\Gvp^{\Ga_k,*} (s)|}\Big| \Gs^{\Ga_k,*}(s)(dv)\\
\nnb &\leq & 2C_\Gve\,\int_U  w^{\Ga_k,*}\ip{ f_y(w^{\Ga_k,*}s,z^{\Ga_k,*}(s),v)  {\Gvp^{\Ga_k,*} (s)\over |\Gvp^{\Ga_k,*} (s)|},{\Gvp^{\Ga_k,*} (s)\over |\Gvp^{\Ga_k,*} (s)|}}  \Gs^{\Ga_k,*}(s)(dv)+2 C_\Gve w^{\Ga_k,*}\\
\nnb &\leq & {2C_\Gve\over c_0^2}\,\int_U  \ip{ w^{\Ga_k,*} f_y(w^{\Ga_k,*}s,z^{\Ga_k,*}(s),v)  \Gvp^{\Ga_k,*} (s), \Gvp^{\Ga_k,*} (s)}  \Gs^{\Ga_k,*}(s)(dv)+ 2w^*C_\Gve
\end{eqnarray*}
when $ d(\Gvp^{\Ga_k,*}(s), Q)< \Gd_\Gve$ and
$$
\Big|\doo {}s {\Gvp^{\Ga_k,*} (s)\over |\Gvp^{\Ga_k,*} (s)|}\Big| \leq  2 w^*\tiN
$$
when $ d(\Gvp^{\Ga_k,*}(s), Q)\geq  \Gd_\Gve$.
Thus, we have
\begin{eqnarray}\label{E371}
\nnb && \Big|\doo {}s {\Gvp^{\Ga_k,*} (s)\over |\Gvp^{\Ga_k,*} (s)|}\Big|  \\
\nnb &\leq & {2C_\Gve\over c_0^2}\,\int_U  \ip{ w^{\Ga_k,*} f_y(w^{\Ga_k,*}s,z^{\Ga_k,*}(s),v)  \Gvp^{\Ga_k,*} (s), \Gvp^{\Ga_k,*} (s)}  \Gs^{\Ga_k,*}(s)(dv)+ 2w^* \tiN\\
 &=&  -{C_\Gve\over c_0^2}\,\doo {} s\Big|\Gvp^{\Ga_k,*} (s)\Big|^2 +2w^* \tiN, \qq \all s\in (s_1,s_2).
\end{eqnarray}
Therefore
\begin{eqnarray}\label{E372}
\nnb && \Big| {\Gvp^{\Ga_k,*} (s_2)\over |\Gvp^{\Ga_k,*} (s_2)|}-{\Gvp^{\Ga_k,*} (s_1)\over |\Gvp^{\Ga_k,*} (s_1)|}\Big|\\
\nnb &\leq &  {C_\Gve\over  c_0^2}\Big(|\Gvp^{\Ga_k,*} (s_1)|^2-|\Gvp^{\Ga_k,*} (s_2)|^2\Big)+2w^*\tiN (s_2-s_1) \\
\nnb &\leq &  {C_\Gve\over  c_0^2}\Big(|\Gvp^{\Ga_k,*} (1-\Gd)|^2-|\Gvp^{\Ga_k,*} (1)|^2\Big)+2w^*\tiN \Gd \\
&\leq & {\Gve\over 2}.
\end{eqnarray}
(b) The condition \refeq{E370} fails but there exists $1\leq j\leq m_\Gve$ such that
\begin{equation}\label{E370b}
{\Gvp^{\Ga_k,*} (s)\over |\Gvp^{\Ga_k,*} (s)|}\not\in \bigcup_{1\leq k\leq m_\Gve \atop k\ne j}E_k, \qq\all s\in (s_1,s_2).
\end{equation}
Then
$$
\ul{s}\equiv \min \set{s\in [s_1,s_2]\Big|{\Gvp^{\Ga_k,*} (s)\over |\Gvp^{\Ga_k,*} (s)|}\in E_j}, \q \bs\equiv \max \set{s\in [s_1,s_2]\Big|{\Gvp^{\Ga_k,*} (s)\over |\Gvp^{\Ga_k,*} (s)|}\in E_j}, \q
$$
are well defined. It holds that
$$
{\Gvp^{\Ga_k,*} (s)\over |\Gvp^{\Ga_k,*} (s)|}\not\in \bigcup_{1\leq k\leq m_\Gve}E_k, \qq\all s\in (s_1,\ul s)\bigcup (\bs, s_2).
$$
By\refeq{E367} and what we get in (a), we have
\begin{eqnarray}\label{E372b}
\nnb && \Big| {\Gvp^{\Ga_k,*} (s_2)\over |\Gvp^{\Ga_k,*} (s_2)|}-{\Gvp^{\Ga_k,*} (s_1)\over |\Gvp^{\Ga_k,*} (s_1)|}\Big|\\
\nnb &\leq &  \Big| {\Gvp^{\Ga_k,*} (s_2)\over |\Gvp^{\Ga_k,*} (s_2)|}-{\Gvp^{\Ga_k,*} (\bs)\over |\Gvp^{\Ga_k,*} (\bs)|}\Big|
+\Big| {\Gvp^{\Ga_k,*} (\bs)\over |\Gvp^{\Ga_k,*} (\bs)|}-{\Gvp^{\Ga_k,*} (\ul s)\over |\Gvp^{\Ga_k,*} (\ul s)|}\Big| +\Big| {\Gvp^{\Ga_k,*} (\ul s)\over |\Gvp^{\Ga_k,*} (\ul s)|}-{\Gvp^{\Ga_k,*} (s_1)\over |\Gvp^{\Ga_k,*} (s_1)|}\Big|\\
&\leq &  {\Gve\over 2}+\Gve+{\Gve\over 2}=2\Gve.
\end{eqnarray}
(c) There exist  $1\leq j_1<j_2\leq m_\Gve$ and $\tis_1,\tis_2\in (s_1,s_2)$ such that
\begin{equation}\label{E370c}
{\Gvp^{\Ga_k,*} (\tis_1)\over |\Gvp^{\Ga_k,*} (\tis_1)|}\in E_{j_1}, \q {\Gvp^{\Ga_k,*} (\tis_2)\over |\Gvp^{\Ga_k,*} (\tis_2)|}\in E_{j_2}.
\end{equation}
Then we can easily see that there must exist $j_3\ne j_4$ and $\tis_3,\tis_4\in (s_1,s_2)$, $\tis_3<\tis_4$ satisfying
\begin{equation}\label{E370c2}
{\Gvp^{\Ga_k,*} (\tis_3)\over |\Gvp^{\Ga_k,*} (\tis_3)|}\in E_{j_3}, \q {\Gvp^{\Ga_k,*} (\tis_4)\over |\Gvp^{\Ga_k,*} (\tis_4)|}\in E_{j_4}
\end{equation}
and
$$
{\Gvp^{\Ga_k,*} (s)\over |\Gvp^{\Ga_k,*} (s)|}\not\in \bigcup_{1\leq k\leq m_\Gve}E_k, \qq\all s\in (\tis_3,\tis_4).
$$
By (a), we have
$$
\Big| {\Gvp^{\Ga_k,*} (\tis_4)\over |\Gvp^{\Ga_k,*} (\tis_4)|}-{\Gvp^{\Ga_k,*} (\tis_3)\over |\Gvp^{\Ga_k,*} (\tis_3)|}\Big|
\leq  {\Gve\over 2}.
$$
Contradicts to \refeq{E368}.  That is, this case should not occur.

Combining the above, we can get that for any $s_1,s_2\in [1-\Gd,1]$,
it always holds that
\begin{equation}\label{EEEEE}
\Big| {\Gvp^{\Ga_k,*} (s_2)\over |\Gvp^{\Ga_k,*} (s_2)|}-{\Gvp^{\Ga_k,*} (s_1)\over |\Gvp^{\Ga_k,*} (s_1)|}\Big|\leq  2\Gve.
\end{equation}
Thus, combining \refeq{EEEEE} with \refeq{E248}, we see that
$\Gvp^{\Ga_k,*}(\cd)$ is equicontinuous on $[0,1]$ and therefore \refeq{Edd1} holds.

Finally, it follows from \refeq{E231a} and \refeq{Edd1} that
\begin{equation}\label{Err}
\ip{ \Gvp^*(1), q-z^*(1)}\geq 0, \qq\all q\in Q.
\end{equation}
That is, \refeq{E357} holds.
\endpf

\begin{Remark}\label{R305} In deriving the transversality condition, the equicontinuity of $z^{\Ga_k,*}(\cd)$ on $[0,1]$ is deducted from $(S3^\pri)$, but have nothing to do with \thb{CE1}---\thb{CE2}. Therefore, $(S3^\pri)$ can be replaced by the following:

$(S3^\prii)$ There exists a homeomorphic mapping $x=x(y)$ such that
$$
g(t,x,u)={\pa x\over \pa y} f_y(t, y(x), u)
$$
is linear increasing in $x$ .
\end{Remark}
\begin{Remark}
The singularity of $f$ near $Q$ bring many difficulties in yielding maximum principle. Though we think the maximum priciple should be hold for any optimal relaxed triple, we failed to got such a result in a general way. Nevertheless, we think it is possible to use the above result to yield
 maximum priciple for an optimal classical control for many special controled systems.
\end{Remark}
\begin{Remark}\label{RRR}  If we know the optimal time $w^*$, then Theorem \ref{T209} still holds if we replce (CE1) by

(CE1$^\pri$) Let $w^*$ be the optimal time of Problem ($\mcR$). There exist   $t_1\in [0,\bt\,]$ and $\Gd_1>0$ such that
$\ds  f_y(t,y,u)+f_y(t,y,u)^\top $
is positive semi-definite for any $t\in [t_1, w^*), d(y,Q)\leq \Gd_1$.

Similarly, (CE2) can be replaced by

 (CE2$^\pri$)
For any $\Gve>0$, there exists $C_\Gve>0$,  $t_\Gve\in (0,w^*)$, $\Gd_\Gve>0$, $m_\Gve\geq 1$  and  $E_1,E_2,\ldots, E_{m_\Gve}\subset S^{n-1}$ such that
\refeq{EEEcde1}---\refeq{EEEcde2} and
and
\begin{equation}\label{EEEbbb}
 |f_y(t,y,u)\xi|\leq C_\Gve\ip{f_y(t,y,u)\xi, \xi}+C_\Gve,  \q\all \xi\in S^{n-1}\setminus \bigcup^m_{k=1} E_k,\q t\in [t_\Gve, w^*), \q d(y,Q)\leq \Gd_\Gve
\end{equation}
for some $C_\Gve>0$ and $t_\Gve\in (0,w^*)$.
\end{Remark}

\def\theequation{4.\arabic{equation}}
\setcounter{equation}{0} 
\setcounter{Definition}{0} \setcounter{Remark}{0}\section{Existence of Optimal Classical Control.}

To get the existence of optimal classical triple, we set the following assumption.

(ES) Assume that for any $(t,y)\in [0,+\infty)\times \IR^n\setminus Q$,
$\ds \set{ f(t,y,u)|u\in U}$ is closed and convex.

\begin{Theorem}\label{T204} Let $y_0\not\in Q$.  Assume \thb{S1}---\thb{S3} and \thb{ES} hold.

If $\mcRP_{ad}\ne \emptyset$, in particular, if $\mcP_{ad}\ne \emptyset$, then Problem \thb{T} admits at least one optimal classical triple. Moreover, any optimal classical triple is also an optimal relaxed triple.
\end{Theorem}
\proof By Theorem \ref{T203}, there exists an optimal relaxed triple $(\bt,\by(\cd),\bGs(\cd))\in\mcRP_{ad}$.
Noting that
$\ds \int_Uf(t, \by(t),u)\bGs(s)(du)$
is the limit of some convex combinations of elements in
$\ds \set{f(t,\by(t),u)|u\in U}$,
we have that for any $t\in [0,\bt)$,
\begin{equation}\label{E21}
\int_Uf(t, \by(t),u)\bGs(s)(du)\in \coh \set{f(t,\by(t),u)|u\in U}.
\end{equation}
Thus it follows from \refeq{E21} and (ES) that
\begin{equation}\label{E21q}
\int_Uf(t, \by(t),u)\bGs(s)(du)\in \set{f(t,\by(t),u)|u\in U},\qq\all t\in [0,\bt).
\end{equation}
Then, by Filippov's lemma, there exists a $\bu(\cd)\in \mcU$ such that
$$
\int_Uf(t, \by(t),u)\bGs(s)(du)=f(t,\by(t),\bu(t)), \qq\eqae t\in [0,\bt\,).
$$
Consequently,
$$
\by(t)=y_0+\int^t_0 f(s,\by(s),\bu(s))\, ds, \qq \all  t\in [0,\bt\,).
$$
This means that $(\bt,\by(\cd),\bu(\cd))\in \mcP_{ad}$ and it is an optimal classical triple.
\endpf

\def\theequation{5.\arabic{equation}}
\setcounter{equation}{0} 
\setcounter{Definition}{0} \setcounter{Remark}{0}\section{Maximum Principle for Optimal Classical Control.}

It is only in the case that an optimal control to Problem (T) is the unique optimal relaxed control to Problem ($\mcR$),
we can get the corresponding maximum principle directly from Theorem \ref{T209}. For many particular systems,
 given an optimal classical triples $(\bt,\by(\cd),\bu(\cd))$, it is possible for us to construct a suitable  $\tif$, such that $(\bt,\by(\cd),\bu(\cd))$ is the unique optimal relaxed triple to a new optimal control problem with
$$
\tif(t,\by(t),\bu(t))= f(t,\by(t),\bu(t)), \q \tif_y(t,by(t),\bu(t))= f_y(t,\by(t),\bu(t)), \qq \eqae t\in [0,\bt\,).
$$
Then, we can get that $(\bt,\by(\cd),\bu(\cd))$ satisfies the maximum principle immediately from Theorem \ref{T209}.

Nevertheless, unfortunately, we fail to prove the maximum principle for every optimal control in general situation.

If the system is affine, namely $f(t,y,u)$ is in the form of $g(t,y)+B(t)u$, then we have the following result.
\begin{Theorem}\label{T210}
Let $y_0\not\in Q$.  Assume $f(t,y,u)=g(t,y)+B(t)u(t)$ and \thb{S1}---\thb{S3} hold, where $B(t)$ takes values in $\IR^{n\times m}$--- the $n\times m$ matrices space. Moreover, assume that $U$ is convex compact subset of $\IR^m$ and  $\mcU_{ad}\ne \emptyset$.
Then Problem \thb{T} admits an optimal classical triple
$(\bt,\by(\cd),\bu(\cd))$ which satisfies
\begin{eqnarray}\label{E231bbb}
&&  \ip{\bpsi(t),   B(t)\bu(t)} =\max_{u\in U }\ip{\bpsi(t),   B(t)u},\q \eqae t\in [0,\bt\,)
\end{eqnarray}
with $\bpsi(\cd)$ being
 a nontrivial solution of
\begin{equation}\label{E227abc}
   \doo {\bpsi(t)} t= \ds -  {\pa g\over \pa y}(t,\by(t))  \bpsi(t),\q    t\in [0,\bt\,).
\end{equation}
\end{Theorem}

\proof We get easily that $\mcU_{ad}\ne \emptyset$ implies $\mcRP_{ad}\ne \emptyset$. By Theorem \ref{T209}, Problem ($\mcR$) admits an optimal relaxed triple
$(\bt,\by(\cd),\bGs(\cd))$ which satisfies
\begin{eqnarray}\label{E254a}
 && \supp \bGs(t)\subseteq \Big\{\tiv\in U\big|  \ip{\bpsi(t), B(t)\tiv}=\max_{v\in U }\ip{\bpsi(t), B(t)v}\Big\},\q \eqae t\in [0,\bt\,)
\end{eqnarray}
with
 $\bpsi(\cd)$ being
 a nontrivial solution of \refeq{E231bbb}.

Since $U$ is convex, by Filippov's lemma, there exists $\bu(\cd)\in \mcU$ such that
\begin{equation}\label{E254b}
B(t)\bu(t)=\int_U B(t)\,v \bGs(t)(dv), \qq\eqae t\in [0,\bt\,).
\end{equation}
Thus
\begin{eqnarray*}
\doo {\by(t)} t &=& \ds \int_U \Big( g(t,\by(t))+B(t)v
\Big)\bGs(t)(dv)\\
&=& \Big(g(t,\by(t))+B(t)\bu(t)\Big), \q \eqae t\in [0,\bt\,).
\end{eqnarray*}
This means $(\bt,\by(\cd),\Gd_{\bu(\cd)})$ is an optimal relaxed triple to Problem ($\mcR$). Equivalently,
$(\bt,\by(\cd),\bu(\cd))$ is an optimal classical triple to Problem \thb{T}.

Finally, \refeq{E231bbb} follows from \refeq{E254a} and \refeq{E254b}.\endpf

\def\theequation{6.\arabic{equation}}
\setcounter{equation}{0} 
\setcounter{Definition}{0} \setcounter{Remark}{0}\section{Optimal time control problems for some particular systems}

In this section, we will show how to apply Theorem \ref{T209} to yield
maximum principle for (every) optimal classical control for particular systems. We will discuss two examples considered in  \cite{Lin1} and \cite{LinWang}.  A crucial property we used here is the `` monotonicity" of the controlled systems.

The first example is concerned optimal quenching time which is considered in \cite{Lin1}.

\noindent \textbf{Example 1.} Let $n=2$, $Q=\set{\pmatrix{x_1\cr x_2}\in \IR^2\Big|x_1=1}$, $U=\set{u\in \IR^2\Big||u|\leq \rho_0}$. Consider the controlled system
\begin{equation}\label{E600}
\left\{\begin{array}{ll}
 \ds \doo{y(t)} t=f(t,y(t),u(t)), &   t\geq 0,\\
 y(0)=y_0 & \end{array}\right.
\end{equation}
with
\begin{equation}\label{E601}
f(t,y,u)=\pmatrix{\ds {y_2\over 1-y_1}\cr y_1+y_2}+B(t)u,   \qq  t\geq 0, y\equiv\pmatrix{y_1\cr y_2}\in \IR^2\setminus Q, u\equiv \pmatrix{u_1\cr u_2}\in U,
\end{equation}
where $y_0=\pmatrix{y_{10}\cr y_{20}}$, $B(\cd)\in L^\infty([0,+\infty);\IR^{2\times 2})$.

\bigskip

We give first a lemma concerning the monotonicity of a system related to \refeq{E601}.

\begin{Lemma}\label{T601} Let $T>0$, $g(\cd)\in L^\infty([0,T]; \IR^2)$. Assume that the solution of equation
\begin{equation}\label{E602}
 \left\{\begin{array}{ll}\ds  \doo {} t \pmatrix{y_1(t)\cr y_2(t)}=\pmatrix{\ds {y_2(t)\over 1-y_1(t)}\cr y_1(t)+y_2(t)}+g(t),&　t>0,\\
   \pmatrix{y_1(0)\cr y_2(0)}=\pmatrix{y_{10}\cr y_{20}}. & \end{array}\right.
\end{equation}
exists on $[0,T)$. Moreover, $\ds y_1(T)\equiv \lim_{t\to T^-}y_1(t)$ exits. Let $h(\cd)\in L^\infty([0,T])$ with $\set{t\in [0,T]|h(t)\ne 0}$ admitting positive measure. Denote $\hy_1(\cd)$ as the solution of equation
\begin{equation}\label{E603}
   \left\{\begin{array}{ll}\ds  \doo {} t\pmatrix{\hy_1(t)\cr \hy_2(t)}=\pmatrix{\ds {\hy_2(t)\over 1-\hy_1(t)}\cr \hy_1(t)+\hy_2(t)}+g(t)+\pmatrix{  h(t)\cr 0},&　t>0,\\
   \pmatrix{\hy_1(0)\cr \hy_2(0)}=\pmatrix{y_{10}\cr y_{20}}. & \end{array}\right.
\end{equation}

\begin{enumerate}\renewcommand{\labelenumi}{\rm(\roman{enumi})}
\item
  If $y_{10}< 1$,
\begin{equation}\label{F604}
h(t)\leq 0, \qq\eqae t\in [0,T],
\end{equation}
then $\hy(\cd)$ exists on $[0,T]$   and $\hy_1(T)<y_1(T)$.

\item
  If $y_{10}> 1$,
\begin{equation}\label{F604}
h(t)\geq 0, \qq\eqae t\in [0,T],
\end{equation}
then $\hy(\cd)$ exists on $[0,T]$   and
$\hy_1(T)> y_1(T)$.
\end{enumerate}
\end{Lemma}
\proof We only prove (i) while (ii) can be proven similarly.
The proof will be finished in two steps.

\textbf{Step I.} Suppose that $y_1(T)\ne 1$. Thus, $y_1(T)<1$ since $y_{10}<1$. Moreover $y(\cd)$ is the solution of \refeq{E603} on $[0,T]$.

Let
$$
\mcV\equiv \set{v(\cd):[0,T]\to [0,1]\Big| v(\cd) \q \mbox{is measurable}}.
$$
Consider the system
\begin{equation}\label{E604}
   \left\{\begin{array}{ll}\ds  \doo {} t\pmatrix{z_1(t)\cr z_2(t)}=\pmatrix{\ds {z_2(t)\over 1-z_1(t)}\cr z_1(t)+z_2(t)}+g(t)+\pmatrix{  h(t)v(t)\cr 0},&　t>0,\\
   \pmatrix{z_1(0)\cr z_2(0)}=\pmatrix{y_{10}\cr y_{20}},  & \end{array}\right.
\end{equation}
and the following optimal control problem: to find $v(\cd)\in  \mcV$ such that $\ds z_1(T)$ takes the minimum value.

Since \refeq{E604}  admits a solution $\ds \pmatrix{y_1(\cd)\cr y_2(\cd)}$ on $[0,T]$ when  $v(\cd)\equiv 0$, we can easily see that the optimal control $\bv(\cd)$ for the above optimal control problem exists by the convexity. We denote the corresponding solution as $\ds \pmatrix{\bz_1(\cd)\cr \bz_2(\cd)}$. By Pontryagin's maximum principle, we have
\begin{equation}\label{E605}
 \bGvp_1(t)  h(t)\bv(t)=\max_{v\in [0,1]}  \Big(\bGvp_1(t) h(t)v \Big) , \qq\eqae t\in [0,T],
\end{equation}
where
\begin{equation}\label{E606}
   \left\{\begin{array}{ll}
   \ds  \doo {\bGvp_1(t)} t =- {\bz_2(t)\over (1-\bz_1(t))^2}\, \bGvp_1(t)- \bGvp_2(t),&　t\in [0,T],\\
   \ds  \doo { \bGvp_2(t)} t =-  {1\over  1-\bz_1(t)} \, \bGvp_1(t)- \bGvp_2(t),&　t\in [0,T],\\
  \bGvp_1(T)=-1, \q  \bGvp_2(T)=0.  & \end{array}\right.
\end{equation}
Denote
$$
\bGP_1(t)=\exp\Big(\int^t_0{\bz_2(s)\over (1-\bz_1(s))^2}\, ds\Big)\bGvp_1(t), \q \bGP_2(t)=e^t\bGvp_2(t), \qq t\in [0,T].
$$
Then \refeq{E606} is converted into
\begin{equation}\label{E607}
   \left\{\begin{array}{ll}
   \ds  \doo {\bGP_1(t)} t =- \exp\Big(\int^t_0{\bz_2(s)\over (1-\bz_1(s))^2}\, ds-t\Big)\bGP_2(t),&　t\in [0,T],\\
   \ds  \doo { \bGP_2(t)} t =-  {1\over  1-\bz_1(t)} \,  \exp\Big(t-\int^t_0{\bz_2(s)\over (1-\bz_1(s))^2}\, ds\Big)\bGP_1(t),&　t\in [0,T],\\
  \bGP_1(T)=-\exp\Big(\int^T_0{\bz_2(s)\over (1-\bz_1(s))^2}\, ds\Big), \q  \bGP_2(T)=0,  & \end{array}\right.
\end{equation}
Since $\bGP_1(T)<0$, there exists a minimum $\eta\in [0,T)$ satisfying
\begin{equation}\label{F608}
\bGP_1(t)<0, \qq \all t\in (\eta,T].
\end{equation}
Consequently,
\begin{equation}\label{F609}
\ds \doo { \bGP_2(t)} t>0, \qq \all t\in (\eta,T].
\end{equation}
Combining with $\bGP_2(T)=0$ we know that
\begin{equation}\label{F610}
\bGP_2(t)<0, \qq t\in (\eta,T].
\end{equation}
Therefore $\bGP_1(\cd)$ is strictly increasing on $[\eta,T]$ and we have
$\bGP_1(\eta)<\bGP(T)<0$. Thus, it must hold that $\eta=0$ and then we get that both $\bGP_1(\cd)$ and $\bGP_2(\cd)$ are negative on $[0,T)$. Thus $\bGvp_1(\cd)$ is negative on $[0,T)$. Finally, it follows from \refeq{E605} that
$$
h(t)\bv(t)=h(t), \qq\eqae t\in [0,T].
$$
This means $\hv(\cd)\equiv 1$ is an optimal control\footnote{We would like to mention that $\bv(\cd)$ is not necessarily $\hv(\cd)$.}, and $\tiv(\cd)\equiv 0$ is not an optimal control. Thus, $\hy_1(T)=\bz_1(T)<y_1(T)$.

\textbf{Step II.} Suppose that $y_1(T)=1$. Denote
$$
Y_1(t)= \hy_1(t)-y_1(t), \q Y_2(t)= \hy_2(t)-y_2(t), \qq t\in [0,T).
$$
Then, by what we get in Step I, we have
\begin{equation}\label{F613A}
Y_1(t)\leq 0, \qq \all t\in (0,T)
\end{equation}
and
\begin{equation}\label{F613}
Y_1(t)<0, \qq \all t\in (S,T),
\end{equation}
where $S\in (0,T)$ such that  $\set{t\in [0,S]|h(t)\ne 0}$ has positive measure.

 We claim that $y_1(T)=1$ implies $y_2(T)\geq 0$.
Since  $\ds \doo { y_2(t)}t$ is bounded, $\ds y_2(T)\equiv \lim_{t\to T^-}y_2(t)$ exists.
If $y_2(T)<0$, then it follows from \refeq{E602} that
$$
\lim_{t\to T^-}\doo{y_1(t)}t=-\infty.
$$
Contradicting the fact of $\ds\lim_{t\to T^-}y_1(t)=1$ and
$$
y_1(t)<1, \qq\all t\in [0,T).
$$
Therefore, we must have $y_2(T)\geq 0$.

On the other hand,
\begin{equation}\label{E613}
 \left\{\begin{array}{ll}\ds \doo {Y_1(t)} t ={y_2(t)\over (1-y_1(t))(1-\hy_1(t)}Y_1(t)+{1\over 1-y_1(t)}Y_2(t)+h(t),& t\in [0,T),   \\
 \doo {Y_2(t)} t =Y_1(t)+Y_2(t),& t\in [0,T),   \\
Y_1(0)=Y_2(0)=0, &
\end{array}\right.
\end{equation}
Combining \refeq{E613} with \refeq{F613A}---\refeq{F613}, we get $Y_2(T)<0$.

If $y_2(T)=0$. Then $\hy_2(T)<0$. This implies $\hy_1(T)<1$.

If $y_2(T)>0$, then there exits $\Gve\in (0,T-S)$ such that
$$
y_2(t)>0, \qq \all t\in (T-\Gve,T).
$$
Therefore
$$
\doo {Y_1(t)} t \leq 0,\qq \all t\in [T-\Gve,T).
$$
Thus $Y_1(T)\leq Y_1(T-\Gve)<0$. That is $\hy_1(T)<1$.

We get the proof.
\endpf

For the time optimal control corresponding to Example 1, we have the following theorem, which is an extension of Theorem 1.4 in \cite{Lin1}.
\begin{Theorem}\label{T602}
Assume that $y_{10}\ne 1$. Let $(\bt,\by(\cd),\bu(\cd))$ be an optimal classical triple of Problem \thb{T}
corresponding to Example 1. Then
\begin{equation}\label{E610b}
 \by_2(\bt\,)\geq 0
\end{equation}
and
there exists a nontrivial solution $\bpsi(\cd)\in C([0,T);\IR^2)$ of the following equation
\begin{equation}\label{E608}
   \doo {\bpsi(t)} t= \ds -\pmatrix{\ds {\by_2(t)\over (1-\by_1(t))^2} &  1\cr \ds
    {1\over  1-\by_1(t)} & 1}  \bpsi(t),\q    t\in [0,\bt\,)
\end{equation}
such that
\begin{eqnarray}\label{E609}
 &&  \ip{\bpsi(t),   B(t)\bu(t)} =\max_{|u|\leq \rho_0 }\ip{\bpsi(t),   B(t)u},\q \eqae t\in [0,\bt\,).
\end{eqnarray}
Moreover,  if $\by_2(\bt\,)>0$, then
\begin{equation}\label{E610}
   \lim_{t\to \bt^-}\bpsi(t)=0.
\end{equation}
\end{Theorem}
\proof We proof only the case of $y_{10}<1$ .

We will get the results by introducing a new system such that the optimal classical triple in consideration becomes the unique optimal relaxed triple of the new optimal problem.

\textbf{I.} Consider that
\begin{equation}\label{E610A}
  \tif(t,y,u)=f(t,y,u)-\pmatrix{\ds  |u-\bu(t)|^2\cr 0},  \q t\geq 0, y\in \IR^2\setminus Q, u\in U.
\end{equation}
We call Problem ($\mcR$) corresponding to $\tif$ as Problem ($\widetilde{\mcR}$).

By direct verification we can see (S1)---(S4) respect to $\tif$ 
hold. 

Let $(t^*,y^*(\cd),\Gs^*(\cd))$ be an optimal relaxed triple of Problem ($\widetilde{\mcR}$). Then $t^*\leq \bt$ since
$(\bt,\by(\cd),\Gd_{\bu(\cd)})$ is also an admissible  relaxed triple of Problem ($\widetilde{\mcR}$).

Let $\hy(\cd)$ be the solution of
$$
\left\{\begin{array}{ll}\ds \doo{\hy(t)}t=  \int_U f (t,\hy(t),u)\Gs^*(t)(du), & t\in [0,t^*), \\
     \hy(0)=y_0, & \end{array}\right.
$$
Then, by Lemma \ref{T601},
$$
\hy_1(t)\geq y_1^*(t), \qq\all t\in [0,t^*).
$$
Since $\ds \lim_{t\to t^*} y^*_1(t)=1$ and
$$
\hy_1(t)<1, \qq\all t\in  [0,t^*)\subseteq (0,\bt\,),
$$
we get
$$
\lim_{t\to t^*} \hy_1(t)=1.
$$
This means
$(t^*,\hy(\cd),\Gs^*(\cd))$ is an admissible relaxed triple of Problem ($\mcR$). Therefore, we have $t^*\geq \bt$ and consequently the equality $t^*=\bt$ holds.

Further, since
$$
\lim_{t\to \bt^-} \hy_1(t)=\lim_{t\to \bt^-} y^*_1(t)=1,
$$
we get from Lemma \ref{T601} that $\set{t\in [0,t_1]| h(t)\ne 0}$ has zero measure,
where
$$
h(t)= \int_U |u-\bu(t)|^2\Gs^*(t)(du)\geq 0, \qq t\in [0,\bt\,).
$$
This means that $\ds \Gs^*(t)=\Gd_{\bu(t)}, \eqae t\in [0,\bt\,)$. Consequently,
\begin{equation}\label{E611}
(y^*(t),\Gs^*(t))=(\by (t),\Gd_{\bu(t)}),  \qq\eqae t\in [0,\bt\,).
\end{equation}

\textbf{II.} We have \refeq{E610b} (see Step II in the proof of Lemma \ref{T601}.

By I, $(\bt, \by (\cd),\Gd_{\bu(\cd)})$ is the unique optimal relaxed triple of ($\widetilde{\mcR}$). Then it follows directly from
Theorem \ref{T209} that there exists a nontrivial solution of \refeq{E608} satisfying \refeq{E609}.

Now, change variable for the system corresponding to $\tif$,
$$
\left\{\begin{array}{l} \ds \tix_1 =(1-y_1 )^2, \\
\tix_2 =y_2.\end{array}\right.
$$
We have
\begin{eqnarray}\label{E601x}
\nnb   \doo{} t \pmatrix{\tix_1(t)\cr \tix_2(t)}   &=& \pmatrix{\ds -2\tix_2(t) \cr 1-\sqrt{\tix_1(t)}+\tix_2(t)}+\pmatrix{-2\sqrt{\tix_1(t)} & 0\cr 0 & 1}\, B(t)u(t) \\
\nnb && \\
\nnb & &  +2\sqrt{\tix_1(t)}\pmatrix{\ds |u(t)-\bu(t)|^2\cr 0}v(t).
\end{eqnarray}
We can verify easily that $(S3^\prii)$ holds.
We can verify easily that $(S3^\prii)$ holds.

If $\by_2(\bt\,)>0$, then there exists $\tit\in (0,\bt\,)$ such that
$$
\by_2(t)>0, \qq\all t\in [\tit,\bt\,].
$$
Then, we can verify that (CE1$^\pri$)---(CE2$^\pri$) hold.

Thus, by Theorem \ref{T209}, Remark \ref{R305} and Remark \ref{RRR},
$\ds \bpsi(\bt\,)=\lim_{t\to \bt^-}\bpsi(t)$ exists and $\ds \bpsi_2(\bt\,)=0$ (see \refeq{E357}).

On the other hand, one can easily see that
$$
\Big|\doo {(1-\by(t))^2} t\Big| \leq C, \qq \all t\in [0,\bt\,]
$$
for some $C>0$. Consequently,
$$
(1-\by(t))^2 \leq C(\bt-t), \qq \all t\in [0,\bt\,].
$$
Then
$$
  \int^{\bt}_0 {\by_2(t)\over (1-\by_1(t))^2}\, dt=+\infty
$$
and it follows from
$$
 \doo {\bpsi_1(t)} t  = -  { \by_2(t)\over (1-\by_1(t))^2}\bpsi_1(t)-\bpsi_2(t)
$$
and L'Hospital's rule that
\begin{eqnarray}
\nnb  \lim_{t\to \bt^-}  \bpsi_1(t) &=&\lim_{t\to \bt^-} {\bpsi_1(0)-\int^t_0 e^{F(s)} \bpsi_2(s)\, ds\over e^{F(t)}}\\
&=& -\lim_{t\to \bt^-} { e^{F(t)} \bpsi_2(t) \over {\by_2(t)\over (1-\by_1(t))^2}e^{F(t)}}=0,
\end{eqnarray}
where
$$
F(t)= \int^t_0 {\by_2(s)\over (1-\by_1(s))^2}\, ds, \qq t\in [0,\bt\,).
$$
That is, we get \refeq{E610}.
\endpf

Now, we consider an example concerning optimal blowup time which is considered in \cite{LinWang}.

\noindent \textbf{Example 2.} Assume $n\geq 1$, $m\geq 1$,  $p>1$ and  $\rho_0>0$. Let $U=\set{u\in \IR^m\Big||u|\leq \rho_0}$. Consider the controlled system
\begin{equation}\label{F624}
\left\{\begin{array}{ll}
 \ds \doo{y(t)} t=|y(t)|^{p-1}y(t)+B(t)u(t), &   t\geq 0,\\
 y(0)=y_0 & \end{array}\right.
\end{equation}
where $B(\cd)\in L^\infty([0,+\infty);\IR^{n\times m})$. Denote by $y(\cd;u(\cd)$ the solution of \refeq{F624} corresponding to $u(\cd)\in \mcU$.
Consider the following blowup time optimal control problem:

\textbf{Problem (BT)}: Find $(\bt,\by(\cd),\bu(\cd))\in \mcP_{ad}$  such that
\begin{equation}\label{F625}
\bt=\inf_{(T,y(\cd),u(\cd))\in \mcP_{ad}}T,
\end{equation}
where
\begin{equation}\label{F626}
\begin{array}{l}\ds
\mcP=\set{(T,y(\cd),u(\cd))\in (0,+\infty)\times C([0,T);\IR^n)\times \mcU\Big| \mbox{\refeq{F624} \,holds on }\, [0,T)},\\
\ds \mcP_{ad}=\set{(T,y(\cd),u(\cd))\in \mcP\Big| \lim_{t\to T^-}|y(t)|=+\infty}.
\end{array}
\end{equation}

\bigskip

Let $(\bt,\by(\cd),\bu(\cd))$ be an optimal triple of Problem (BT).
To yield the maximum principle, we consider the following modified relaxed system:
\begin{equation}\label{F627}
\left\{\begin{array}{ll}
 \ds \doo{y(t)} t=|y(t)|^{p-1}y(t)+ \int_U \Big(B(t)u-{|u-\bu(t)|^2\over 4\rho_0^2} y(t)\Big)\Gs(t)(du), &   t\geq 0,\\
 y(0)=y_0. & \end{array}\right.
\end{equation}

We give some observations first. Denote
\begin{equation}\label{F627B}
M=\sup_{(t,u)\in [0+\infty)\times U}|B(t)u|.
\end{equation}
For  $s>0$, $r>0$, let
$\ol\GT(\cd;s,r)$ and $\ul\GT(\cd;s,r)$ be the solution of
\begin{equation}\label{F627}
\left\{\begin{array}{l}\ds\doo {\Gt(t)} t =\Gt^p(t)+M+\Gt(t),   \\
\Gt(s)=r,
\end{array}\right.
\end{equation}
and
\begin{equation}\label{F628}
\left\{\begin{array}{ll}\ds\doo {\Gt(t)} t =\Gt^p(t)-M-\Gt(t),   \\
\Gt(s)=r
\end{array}\right.
\end{equation}
respectively.  Let $r_0>1$ be big enough, say for example, let
\begin{equation}\label{F628B}
r_0={p+M\over p-1}.
\end{equation}
Then it is easy to see that for any $r>r_0$,
$\ol\GT(\cd)$ and $\ul\GT(\cd)$ will be blowup at $t=\Xi^*(r)+s$ and $t=\Xi_*(r)+s$, respectively.
Then for any $y(\cd)$ satisfies
\begin{equation}\label{F629B}
 \ds \doo{y(t)} t=|y(t)|^{p-1}y(t)+ \int_U \Big(B(t)u-h(t) y(t)\Big)\Gs(t)(du)
\end{equation}
with $r\equiv |y(s)|>r_0$ and
$$
\|h(t)\|\leq 1, \qq\all t>0,
$$
it holds that
\begin{equation}\label{F630}
\ul\GT(t)<|y(t)|< \ol\GT(t;s,r), \qq  t>s,
\end{equation}
and
\begin{equation}\label{F631}
\ul\GT(t;s,r)>|y(t)|> \ol\GT(t;s,r), \qq  s_0<t<s,
\end{equation}
where
$$
s_0=\inf\set{t>0|\ol\GT(t;s,r)>r_0}.
$$
Since both $\Xi^*(\cd)$ and $\Xi_*(\cd)$ are strictly decreasing on $(r_0,+\infty)$,
they have inverse functions $\xi^*(\cd)$ and $\xi_*(\cd)$.
In fact, we have
\begin{equation}\label{F633B}
\left\{\begin{array}{ll}\ds
\xi^*(\Xi^*(r)-t)=\ol\GT(t;s,r), & t\in (0,\Xi^*(r)-s), \\
\ds \xi_*(\Xi^*(r)-t)=\ul\GT(t;s,r), & t\in (0,\Xi_*(r)-s),\end{array}\right.  \qq \all r>r_0.
\end{equation}

We establish a monotonicity lemma related to \refeq{F629B}.
\begin{Lemma}\label{T602} Assume that $T>s>0$, $p-1>\Ga>0$, $g(\cd)\in L^\infty([s,T]; \IR^n)$,   $h(\cd)\in L^\infty([s,T])$,
\begin{equation}\label{F635B}
|y_s|\geq\tiM\equiv \max\Big[r_0, \Big({4\Ga\over p-1-\Ga}\Big)^{1\over p-1}, \Big({2e^{\Xi(r_0)}\over p-1-\Ga}\Big)^{1\over p}, 2e^{\Xi_*(r_0)}M\Big],
\end{equation}
where $M$ and $r_0$  are defined by \refeq{F627B} and \refeq{F628B},
\begin{equation}\label{F635}
|g(t)|\leq M, \q 0\leq h(t)\leq 1, \qq t\in [s,T]
\end{equation}
and $\set{t\in [s,T]|h(t)\ne 0}$ admits positive measure.

Suppose that $\tiy(\cd)$ is the solution of equation
\begin{equation}\label{F636}
\left\{\begin{array}{ll}
 \ds \doo{y(t)} t=|y(t)|^{p-1}y(t)+g(t), &   t\geq s,\\
 y(s)=y_s, & \end{array}\right.
\end{equation}
which exists on $[s,T]$.  Denote by $\hy(\cd)$ the solution of equation
\begin{equation}\label{F639}
\left\{\begin{array}{ll}
 \ds \doo{y(t)} t=|y(t)|^{p-1}y(t)+g(t)-h(t)y(t), &   t\geq s,\\
 y(s)=y_s & \end{array}\right.
\end{equation}
then $\hy(\cd)$ exists on $[s,T]$. Moreover,
\begin{equation}\label{F640}
|\hy(T)|^{-\Ga}\geq \int^T_s \Ga |\xi_*(T-t)|^{-\Ga} h(t)\, dt.
\end{equation}
\end{Lemma}
\proof We will prove the lemma in two steps.

\textbf{I.}
Let
$$
\mcV\equiv \set{v(\cd):[s,T]\to [0,1]\Big| v(\cd) \q \mbox{is measurable}}.
$$
Consider the system
\begin{equation}\label{F641}
\left\{\begin{array}{ll}
 \ds \doo{y(t)} t=|y(t)|^{p-1}y(t)+g(t)-h(t)y(t)v(t), &   t\geq s,\\
 y(s)=y_s. & \end{array}\right.
\end{equation}
For $v(\cd)\in \mcV$, denote $y(\cd;v(\cd))$ the corresponding solution of \refeq{F641}.
It is easy to see that \refeq{F635B} implies
\begin{equation}\label{FFFF}
|y(t;v(\cd)|\geq \tiM, \qq\all t\in [s,S], v(\cd)\in \mcV.
\end{equation}
We will prove that if for some $S\in (0,T]$,
\begin{equation}\label{F642}
\sup_{v(\cd)\in \mcV\atop t\in [s,S]}|y(t;v(\cd))|<+\infty,
\end{equation}
then
\begin{equation}\label{F643}
|y(S;v(\cd))|^{-\Ga}\geq |\tiy(S)|^{-\Ga}+\int^S_s \Ga |\xi_*(S-t)|^{-\Ga} h(t)v(t)\, dt, \qq \all v(\cd)\in \mcV.
\end{equation}
Denote
\begin{equation}\label{F643B}
\xi(t)=\Ga |\xi_*(S-t)|^{-\Ga}, \qq t\in [s,S].
\end{equation}
By classical optimal control theory, there exists an optimal control $\bv(\cd)\in \mcV$ such that
\begin{eqnarray}\label{F644}
 \nnb && |y(S;\bv(\cd))|^{-\Ga}-\int^S_s \xi(t) h(t)\bv(t)\, dt\\
 &\leq & |y(S;v(\cd))|^{-\Ga}-\int^S_s \xi(t) h(t)v(t)\, dt, \qq \all v(\cd)\in \mcV.
\end{eqnarray}
Moreover, by Pontryagin's maximum principle, we have
\begin{equation}\label{F645}
h(t)\Big(\xi(t)-\ip{\bGvp (t), \by(t)}\Big)\bv(t) =\max_{v\in [0,1]}  h(t)\Big(\xi(t)-\ip{\bGvp (t), \by(t)}\Big)v , \qq\eqae t\in [s,S],
\end{equation}
where $\by(\cd)=y(\cd;\bv(\cd))$ and
\begin{equation}\label{F646}
   \left\{\begin{array}{ll}
   \ds  \doo {\bGvp(t)} t =- |\by(t)|^{p-1}\Big(I+(p-1){\by(t)\by(t)^\top\over |\by(t)|^2}\Big)\bGvp(t)+h(t)\bv(t)\bGvp(t),&　t\in [s,T],\\
  \bGvp(S)=\Ga |\by(S)|^{-\Ga-2}\by(S).  & \end{array}\right.
\end{equation}
We have
\begin{equation}\label{F647}
   \left\{\begin{array}{ll}
   \ds  \doo{} t \ip{\bGvp(t),\by(t)}  =-(p-1) |\by(t)|^{p-1} \ip{\bGvp(t),\by(t)}+\ip{\bGvp(t),g(t)},&　t\in [s,T],\\
  \ip{\bGvp(S),\by(S)}=\Ga |\by(S)|^{-\Ga}.  & \end{array}\right.
\end{equation}
Let $s_0\in [s,S)$ be the smallest number such that $\ip{\bGvp(\cd),\by(\cd)}$ is non-negative on $[s_0,S]$.
Then
\begin{eqnarray}\label{F648}
\nnb \doo {} t \Big(|\bGvp(t)|\,|\by(t)| \Big)&=& |\bGvp(t)| \ip{{\by(t)\over |\by(t)|}, |\by(t)|^{p-1}\by(t)+g(t)-h(t)\by(t)\bv(t)}\\
\nnb && + |\by(t)| \ip{{\bGvp(t)\over |\bGvp(t)|}, - |\by(t)|^{p-1}\Big(I+(p-1){\by(t)\by(t)^\top\over |\by(t)|^2}\Big)\bGvp(t)+h(t)\bv(t)\bGvp(t)}\\
\nnb &=& |\bGvp(t)| \ip{{\by(t)\over |\by(t)|}, g(t)} -(p-1) |\by(t)|^{p-2} {\ip{\bGvp(t),\by(t)}^2\over |\bGvp(t)|}\\
\nnb &\geq & |\bGvp(t)| \ip{{\by(t)\over |\by(t)|}, g(t)} -(p-1) e^{(S-t)}|\by(t)|^{p-1} \ip{\bGvp(t),\by(t)}\\
\nnb &= & |\bGvp(t)| \ip{{\by(t)\over |\by(t)|}, g(t)} +e^{(S-t)}\Big[\doo{} t \ip{\bGvp(t),\by(t)}-\ip{\bGvp(t),g(t)}\Big]\\
 &\geq & e^{(S-t)}\doo{} t \ip{\bGvp(t),\by(t)}- |\bGvp(t)|\,\by(t)|, \qq t\in [s_0,S].
\end{eqnarray}
Therefore,
\begin{equation}\label{F649}
\doo {} t \Big(e^{t-S}|\bGvp(t)|\,|\by(t)| \Big)\geq   \doo{} t \ip{\bGvp(t),\by(t)},  \qq t\in [s_0,S].
\end{equation}
We get
\begin{equation}\label{F650}
    |\bGvp(S)|\,|\by(S)|   - e^{(t-S)}  |\bGvp(t)|\,|\by(t)| \geq     \ip{\bGvp(S),\by(S)}-  \ip{\bGvp(t),\by(t)}, \qq t\in [s_0,S],
\end{equation}
Thus, it follows from \refeq{F647} and \refeq{F650} that
\begin{equation}\label{F652}
     |\bGvp(t)|\,|\by(t)| \leq e^{\Xi_*(r_0)}  \ip{\bGvp(t),\by(t)}, \q t\in [s_0,S].
\end{equation}
Now, let $\Gz(\cd)=\Ga |\by(\cd)|^{-\Ga}-\ip{\bGvp(\cd),\by(\cd)}$. Then $\Gz(S)=0$. By  \refeq{FFFF} and \refeq{F652}, there is
\begin{eqnarray}
\nnb \doo {\Gz(t)} t &=& -\Ga^2 |\by(t)|^{-\Ga-1}\ip{{\by(t)\over |\by(t)|}, |\by(t)|^{p-1}\by(t)+g(t)-h(t)\by(t)\bv(t)}\\
\nnb &  & +(p-1) |\by(t)|^{p-1} \ip{\bGvp(t),\by(t)}-\ip{\bGvp(t),g(t)}\\
\nnb &\geq & -\Ga^2 |\by(t)|^{-\Ga-1}\Big(|\by(t)|^p+M+|\by(t)|\Big)\\
\nnb && +(p-1) |\by(t)|^{p-1} \ip{\bGvp(t),\by(t)}-M|\bGvp(t)|\\
\nnb &\geq & -{\Ga (\Ga+p-1)\over 2} |\by(t)|^{-\Ga-1} |\by(t)|^p   +{\Ga +p-1\over 2} |\by(t)|^{p-1} \ip{\bGvp(t),\by(t)}\\
&=& -{\Ga +p-1\over 2}\Gz(t), \qq\qq t\in [s_0,S].
\end{eqnarray}
Consequently,
\begin{equation}\label{F655}
  \ip{\bGvp(t),\by(t)}\geq \Ga |\by(t)|^{-\Ga},  \qq t\in [s_0,S].
\end{equation}
Thus, it must hold that $s_0=s$. Moreover, we get from \refeq{F655} and \refeq{F635B} that
\begin{equation}\label{F656}
  \ip{\bGvp(t),\by(t)}\geq \Ga |\xi_*(\Xi_*(|y_s|)-t)|^{-\Ga}>\Ga |\xi_*(S-t)|^{-\Ga}=\xi(t),  \qq t\in [s,S].
\end{equation}
Thus \refeq{F645} implies
\begin{equation}\label{F657}
\bv(t)=0,\qq \eqae t\in \set{\tau\in [s,S]\Big| h(\tau)\ne 0}.
\end{equation}
Consequently, $\by(\cd)=\tiy(\cd)$ and we get \refeq{F643}.

\textbf{II.} Let
\begin{equation}\label{F658}
E=\set{S\in [s,T]\Big|\, \mbox{\refeq{F642} holds}}.
\end{equation}
Then obviously, $E$ is closed and $[s,s+\Xi^*(|y_s|))\bigcap [s,T]\subseteq E$.
On the other hand, if $S\in E$,   we have \refeq{F643}. Consequently,
\begin{equation}\label{F669}
|y(S;v(\cd))|\leq R\equiv \Big[|\tiy(S)|^{-\Ga}+\int^S_s \Ga |\xi_*(S-t)|^{-\Ga} h(t)v(t)\, dt\Big]^{-{1\over \Ga}}, \qq \all v(\cd)\in \mcV.
\end{equation}
This implies $[S,S+\Xi^*(R))\bigcap [s,T]\subseteq E$.  Therefore,  $E$ is also a relatively open set respect to $[s,T]$. We must have $E=[s,T]$.

Replacing $S$ by $T$ in \refeq{F643} we get
\begin{equation}\label{F670}
|y(T;v(\cd))|^{-\Ga}\geq |\tiy(T)|^{-\Ga}+\int^T_s \Ga |\xi_*(T-t)|^{-\Ga} h(t)v(t)\, dt, \qq \all v(\cd)\in \mcV.
\end{equation}
Taking $v(\cd)\equiv 1$, we get \refeq{F640}.
\endpf

Now, we establish maximum principle for solutions of Problem (BT).
\begin{Theorem}\label{T606}
Assume that   $(\bt,\by(\cd),\bu(\cd))$ is an optimal classical triple of Problem \thb{BT}.
Then, there exists a nontrivial solution  $\bpsi(\cd)\in C([0,\bt\,);\IR^n)$ of the following equation
\begin{equation}\label{G319}
\doo {\bpsi(t)} t=-|\by(t)|^{p-1}\Big(I+(p-1){\by(t)\by(t)^\top\over |\by(t)|^2}\Big)\bpsi(t),\q    t\in [0,\bt\,)
\end{equation}
such that
\begin{eqnarray}\label{G320}
 &&  \ip{\bpsi(t),   B(t)\bu(t)} =\max_{ u\in U}\ip{\bpsi(t),  B(t)u},\q \eqae t\in [0,\bt\,).
\end{eqnarray}
Moreover,
\begin{equation}\label{G321}
   \lim_{t\to \bt^-}\bpsi(t)=0.
\end{equation}
\end{Theorem}
\proof Let $\Ga\in (0,p-1)$ and $\tiM$ be defined by
\refeq{F635B} and
\begin{equation}\label{G321B}
  s=\max(0,\Xi^*(\tiM)).
\end{equation}

Consider the following system
\begin{equation}\label{F667}
\left\{\begin{array}{ll}
 \ds \doo{y(t)} t=|y(t)|^{p-1}y(t)+B(t)u(t)-h(t)\chi_{[s,+\infty)]}(t)y(t), &   t>0,\vspace{2mm}\\
 \ds \doo{h(t)} t={|u-\bu(t)|\over \bt\, (1+|u-\bu(t)|)}, & t>0,\\
 y(0)=y_0, h(0)=0,  & \end{array}\right.
\end{equation}
Let $(t^*,(y^*(\cd),h^*(\cd)), \Gs^*(\cd))$ be an optimal triple that makes the solution of \refeq{F667} blowup in its shortest time.
Naturally, we have $t^*\leq \bt$. Then it must hold that
\begin{equation}\label{G321C}
\tiM\leq |y^*(s)|<+\infty.
\end{equation}
We claim that
\begin{equation}\label{F668}
\Gs^*(t)=\Gd_{\bu(t)}, \qq \eqae t\in [0,t^*).
\end{equation}
Otherwise,
$$
1\leq h^*(t)\leq 1, \qq t\in [0,t^*)
$$
and $\set{t\in [0,t^*)| h^*(t)>0}$ has positive measure.  Noting that
\begin{equation}\label{F667BB}
 \doo{y^*(t)} t=|y^*(t)|^{p-1}y^*(t)+B(t)\int_Uu\Gs^*(t)(du)-h^*(t) y^*(t), \qq    t\in (s,t^*),
\end{equation}
we get from Lemma \ref{T602} that
\begin{equation}\label{F667BB}
 |y^*(T)|^{-\Ga}\geq \int^T_s \Ga |\xi_*(T-t)|^{-\Ga} h^*(t)\, dt, \qq T\in (s,t^*).
\end{equation}
This implies
\begin{equation}\label{F672}
\limsup_{T\to t^*} |y^*(T)|<+\infty.
\end{equation}
Contradicts to $y^*(\cd)$ blowup at $t^*$.  Therefore, \refeq{F668} holds and consequently,
\begin{equation}\label{F673}
t^*=\bt, \q y^*(\cd)=\by(\cd).
\end{equation}
This means $(\bt,(\by(\cd),0),\bu(\cd))$ is the unique optimal relaxed triple for shortest blowup optimal time problem respect to system \refeq{F667}.

Now, let $X_0\in \IR^n$ and $\Gb>0$ satisfy
$$
|\by(t)-X_0|>\Gb, \qq |t\in [0,\bt\,).
$$
We can construct a smooth  inversable  map $G:\IR_n\setminus B_\Gb(X_0)\to B_R(0)$ such that
$$
G(y)=|y|^{-\Gg-1}y(t), \qq |y|\geq R_1,
$$
where $R>R_1>Gb$, $\Gg>p-1$, and $B_r(X)$ denotes the ball of radius $r$ and centered at $X$.
Denote $g$ the inverse map of $G$.
Changing variable $z=G(y)$, or equivalently, $y=g(z)$, \refeq{F667} becomes
\begin{equation}\label{F674}
\left\{\begin{array}{l}
 \ds \doo{z(t)} t=\tif(t,z(t),u(t))-h(t)\chi_{[s,+\infty)}(t)\Big({\pa G\over \pa z}(z(t))\Big)^\top  g(z(t)),  \vspace{2mm}\\
 \ds \doo{h(t)} t={|u-\bu(t)|\over \bt\, (1+|u-\bu(t)|)},\\
 z(0)=G(y_0), h(0)=0,    \end{array}\right.
\end{equation}
where
\begin{equation}\label{F675}
\tif(t,z,h,u)=\Big({\pa G(z)\over \pa z}\Big)^\top \Big[|g(z)|^{p-1}g(z)+B(t)u\Big].
\end{equation}
Let $\bz(\cd)=G(\by(\cd))$, $Q=\set{(z,h)|z=0}$. Then $(\bt,(\bz(\cd),0),\bu(\cd))$ becomes the unique optimal relaxed triple for Problem (T) respect to system \refeq{F675}. One can verify that the corresponding assumptions (S1)---(S4) hold. Then, Theorem \ref{T209} implies
\begin{equation}\label{F676}
\ip{\Gvp^*(t), \tif(t,\bz(t),\bu(t))}=\max \ip{\Gvp^*(t), \tif(t,\bz(t),u)}, \qq eqae t\in [0,\bt\,)
\end{equation}
for some nontrivial solution $\Gvp^*(\cd)$ of
\begin{equation}\label{F676}
\doo {\Gvp^*(t)} t= -{\pa  \tif\over \pa z}(t,\bz(t),\bu(t)), \qq t\in [0,\bt\,).
\end{equation}
Let
\begin{equation}\label{F677}
\bpsi(t)={\pa G(y)\over \pa y} \Gvp^*(t)\Big|_{y=\by(t)},  \qq t\in [0,\bt\,).
\end{equation}
Then $\bpsi(\cd)$ is nontrivial and we get \refeq{G319}---\refeq{G320} from \refeq{F677} and \refeq{F676}.

By \refeq{G319},
\begin{eqnarray}\label{E679}
\nnb  {1\over 2}\doo{|\bpsi(t)|^2} t & =& -|\by(t)|^{p-1} |\bpsi(t)|^2-(p-1)|\by(t)|^{p-3}\ip{\by(t),\bpsi(t)}^2\\
 &\leq & -|\by(t)|^{p-1} |\bpsi(t)|^2,\q    t\in [T_0,\bt\,).
\end{eqnarray}
On the other hand, since
\begin{eqnarray}\label{E680}
\nnb  && {1\over 2}\doo{ |\by(t)|^2} t = |\by(t)|^{p-1}|\by(t)|^2+\ip{\by(t),B(t)\bu(t)}\\
  &\leq &    |\by(t)|^{p-1}|\by(t)|^2 +M|\by(t)|, \qq\all t\in [0,\bt\,)
\end{eqnarray}
and $\ds \lim_{t\to \bt^-}|\by(t)|=+\infty$, we should have
\begin{equation}\label{E681}
\lim_{t\to \bt^-}\int^t_0 |\by(\tau)|^{p-1}\, d\tau=+\infty.
\end{equation}
Therefore, \refeq{G320} follows from \refeq{E679} and \refeq{E681}.
\endpf

\end{document}